\documentstyle[amscd,amssymb,verbatim]{amsart}
\newcommand{\lrar}[1]{\begin{picture}(50,10)(-25,-5)                          
\put(-25,0){\vector(1,0){50}}
\put(0,5){\makebox(0,0)[b]{\mbox{$#1$}}}
\end{picture}}

\newcommand{\ldar}[1]{\begin{picture}(10,50)(-5,-25)
\put(0,25){\vector(0,-1){50}}
\put(5,0){\mbox{$#1$}}
\end{picture}}

\newcommand{\Nm}{\operatorname{Nm}}
\newcommand{\veps}{\varepsilon}

\newcommand{\DD}{{\cal D}}
\newcommand{\KK}{{\cal K}}

\newcommand{\SL}{\operatorname{SL}}
\newcommand{\G}{{\Bbb G}}

\newcommand{\gog}{{\frak g}}

\newcommand{\lan}{\langle}
\newcommand{\ran}{\rangle}

\newcommand{\CC}{{\cal C}}

\newcommand{\Tr}{\operatorname{Tr}}
\newcommand{\HS}{{\frak H}}

\newcommand{\Th}{\Theta}
\newcommand{\si}{\sigma}

\newcommand{\eps}{\epsilon}
\newcommand{\We}{\bigwedge}
\renewcommand{\ker}{\operatorname{ker}}

\numberwithin{equation}{section}

\newtheorem{thm}{Theorem}
\newtheorem{prop}{Proposition}[section]
\newtheorem{lem}[prop]{Lemma}
\newtheorem{corint}[thm]{Corollary}
\newtheorem{cor}[prop]{Corollary}
\newenvironment{rem}{\vspace{3mm}\noindent
{\bf Remark.}}{\vspace{3mm}}
\newenvironment{rems}{\vspace{3mm}
\noindent {\bf Remarks.}}{\vspace{3mm}}

\newcommand{\Pf}{\noindent {\it Proof}}
\newcommand{\id}{\operatorname{id}}

\newcommand{\ov}{\overline}

\renewcommand{\Im}{\operatorname{Im}}

\newcommand{\rk}{\operatorname{rk}}
\newcommand{\ra}{\rightarrow}

\newcommand{\FF}{{\cal F}}

\newcommand{\PP}{{\cal P}}

\newcommand{\SS}{{\cal S}}
\newcommand{\LL}{{\cal L}}
\renewcommand{\O}{{\cal O}}
\newcommand{\Om}{\Omega}
\newcommand{\dbar}{\overline{\partial}}
\newcommand{\Hom}{\operatorname{Hom}}
\newcommand{\Ext}{\operatorname{Ext}}

\newcommand{\Res}{\operatorname{Res}}

\renewcommand{\a}{\alpha}
\renewcommand{\b}{\beta}
\newcommand{\om}{\omega}
\newcommand{\De}{\Delta}
\newcommand{\la}{\lambda}
\newcommand{\th}{\theta}
\newcommand{\C}{{\Bbb C}}
\newcommand{\R}{{\Bbb R}}
\newcommand{\Z}{{\Bbb Z}}
\newcommand{\Q}{{\Bbb Q}}
\newcommand{\nn}{{\bf n}}
\newcommand{\cc}{{\bf c}}
\newcommand{\zz}{{\bf z}}
\newcommand{\vv}{{\bf v}}
\newcommand{\ww}{{\bf w}}
\newcommand{\xx}{{\bf x}}
\newcommand{\yy}{{\bf y}}
\newcommand{\uu}{{\bf u}}
\newcommand{\La}{\Lambda}
\newcommand{\Ga}{\Gamma}

\newcommand{\wt}{\widetilde}

\newcommand{\sign}{\operatorname{sign}}
\newcommand{\sub}{\subset}
\newcommand{\ed}{\qed\vspace{3mm}}

\title{Indefinite theta series of signature $(1,1)$ from the point
of view of homological mirror symmetry}
\author{A. Polishchuk}
\address{Department of Mathematics, University of Oregon, Eugene, OR 97403}
\email{apolish@@math.uoregon.edu}

\begin{document}

\begin{abstract} We apply the homological mirror symmetry
for elliptic curves to the study of indefinite theta series. We prove
that every such series corresponding to a quadratic form of signature (1,1)
can be expressed in terms of theta series associated with split
quadratic forms and the usual theta series. We also show that indefinite
theta series corresponding to univalued Massey products between line
bundles on elliptic curve are modular.
\end{abstract}

\maketitle

\bigskip

\centerline{\sc Introduction}

\bigskip

The classical theta series are the series of the form
$\sum_{\nn\in\La}\exp(\pi i\tau Q(\nn)+2\pi i \nn\cdot \zz)$
where $Q$ is a positive definite integer-valued quadratic form
on a lattice $\La$, $\zz\in\La_{\C}$, 
$\zz\cdot \zz':=\frac{1}{2}(Q(\zz+\zz')-Q(\zz)-Q(\zz'))$
is the symmetric
pairing on $\La_{\C}$ induced by $Q$, $\tau$ belongs to the upper half-plane
$\HS$. It is well-known that they are Jacobi forms on $\C\times\HS$ of
weight $\rk\La/2$ (see \cite{EZ}).
Now assume that $\La$ is a rank 2 lattice equipped with a non-degenerate
$\Q$-valued quadratic form $Q$ of signature $(1,1)$.
Let us fix an open cone $C\subset\La_{\R}$ of the form
$C=\{\vv\in\La_{\R}:\phi(\vv)\cdot\psi(\vv)>0\}$ for a pair of linear forms
$\phi$ and $\psi$ on $\La_{\R}$ defined over $\Q$, such that $Q|_C>0$.  Let
$C=C^+\cup C^-$ be the decomposition of $C$ into two connected components,
$\sign:C\ra\{\pm 1\}$ be the corresponding sign function (which is equal to
$1$ on $C^+$ and to $-1$ on $C^-$).
Also let $\a:\La_{\C}\ra\La_{\R}$ be the map given by
$\a(\zz)=\Im(\zz)/\Im(\tau)$.
Then we define the {\it indefinite theta series}\footnote{
The definition of indefinite theta series by G\"ottsche and Zagier
in \cite{GZ} differs slightly from ours in that they fix a connected component
of the domain of definition of $\Th_{\La,Q,C}$. Also we allow
$Q$ to take rational values; however, rescaling $\tau$ and $z$ one
can always reduce to the case of integer-valued form $Q$.} 
associated with $(\La,Q,C)$
by the formula
\begin{equation}\label{maindef}
\Theta_{\La,Q,C}(\zz,\tau)= 
\sum_{\nn\in\La:\nn+\a(\zz)\in C}\sign(\nn+\a(\zz))\exp(\pi i\tau Q(\nn)+
2\pi i \nn\cdot \zz).
\end{equation}
This is a holomorphic function of $\tau$ (in the upper half-plane) and
of the second variable $\zz\in\La_{\C}$ which varies
in the complement to $\a^{-1}(\partial C+\La)$
where $\partial C$ is the boundary of $C$.

Our interest in the functions $\Theta_{\La,Q,C}(\zz,\tau)$ is motivated
by the observation that for some special choices of $\La$, $Q$, $C$,
vectors $\vv,\ww\in\La_{\Q}$ and a rational number $\la$ the function
$$\exp(\pi i\la\tau)\Theta_{\La,Q,C}(\tau\vv+\ww,\tau)$$
is a modular form of weight $1$ for some congruenz-subgroup of $\SL_2(\Z)$.
In particular, the indefinite theta series defined by Hecke in \cite{He}
and their generalizations to arbitrary integral lattices of signature
$(1,1)$ can be written in this form. By definition these series are
$$\Th^H_{\La,Q;\cc}(\tau)=
\sum_{\nn\in \La+\cc/G; Q(\nn)>0}\sign(\nn)\exp(\pi i\tau Q(\nn)),$$
where we assume that $Q(\nn)\in 2\Z$ for all $\nn\in\La$,  
$\cc\in \La^{\perp}=\{\xx\in\La_{\Q}:\xx\cdot\La\subset\Z\}$,
the function $\sign(\nn)$ takes opposite values $\pm 1$ on two
components of the cone $Q>0$, $G$ is the subgroup in the 
identity component of the automorphism group of $Q$ 
consisting of elements preserving $\La+\cc$.
It is easy to see that such a series can be rewritten in the above
form for some rational cone $C$: for this one can use the decomposition
of the set of lattice points in the cone $Q>0$ described in \cite{Z1}.

Another example of modular behaviour of indefinite theta series
goes back to Kronecker: one should consider the quadratic form
$(m,n)\mapsto mn$ on $\Z^2$ and the cone
$C=\{(x,y): xy>0\}$ (see \cite{W},\cite{Z2}). The corresponding
series is a meromorphic Jacobi form on $\C\times\C\times\HS$. 
This example is related to ``Teilwerte'' of Weierstrass zeta-function
considered by Hecke in \cite{He}. 

The main point we would like to make in this paper is
that these examples provide an evidence for the following conjecture:
{\it an indefinite theta series of signature $(1,1)$ is modular
if and only if it corresponds to a universal univalued triple Massey product
in the derived categories of coherent sheaves on elliptic curves}.
This correspondence which is based on homological mirror symmetry for
elliptic curves (proved for transversal products in \cite{P-hmc})
leads to explicit rational expressions for modular indefinite theta
series in terms of the usual theta functions. 

When the form $Q$ is a
product of rational linear forms, the function
$\Theta_{\La,Q,C}(\zz,\tau)$
(restricted to a connected component of its domain of definition)
can be expressed via the following bilateral basic
hypergeometric series
\begin{equation}\label{kappa}
\kappa(y,x;\tau)=\sum_{m\in\Z}\frac{\exp(\pi i\tau m^2 +2\pi i mx)}
{\exp(2\pi i m\tau)-\exp(2\pi i y)}.
\end{equation}
This series was introduced by M.~P.~Appell
in his work \cite{Ap} on decomposition of elliptic functions of the
third kind into simple elements (see also \cite{P-ap}).

As the first application of our techniques
we show that all indefinite theta series associated
with quadratic forms of signature (1,1) and rational cones
can be expressed in terms of $\kappa$ and the usual theta series.
To formulate this more precisely let us call a meromorphic function 
$\phi(z,\tau)$ on $\C\times\HS$ {\it elliptic} if it can be expressed
rationally over $\C$ in terms of functions of the form
$\theta_c(az+b\tau,d\tau)$ and $\exp(\pi i d\tau)$, where 
$a,b,c\in\Q$, $d\in\Q_{>0}$, $\theta_c(z,\tau)=\sum_{n\in\Z+c}\exp
(\pi i\tau n^2+2\pi i nz)$. Let us also denote
$$\kappa_c(y,x;\tau)=\sum_{m\in\Z+c}\frac{\exp(\pi i\tau m^2 +2\pi i mx)}
{\exp(2\pi i m\tau)-\exp(2\pi i y)},$$
where $c\in\Q$.

\begin{thm}\label{main}
For every triple $(\La,Q,C)$ as above and every connected open subset 
$U\subset\La_{\C}\setminus \a^{-1}(\partial C+\La)$ there exist
$\Q$-linear functionals $(r,s;l_i,i=1,\ldots,N)$
on $\La_{\Q}$, constants 
$(a_i,b_i,c_i,d_i,e_i,i=1,\ldots,N)$ in $\Q$, $f\in\Q_{>0}$, 
and meromorphic elliptic functions $(\phi_i,\psi_i,i=1,\ldots,N)$,
such that
$$\Theta_{\La,Q,C}(\zz,\tau)= \sum_{i=1}^N \phi_i(r(\zz))\psi_i(s(\zz))
\kappa_{e_i}(a_i\tau+b_i,l_i(\zz)+c_i\tau+d_i;f\tau)$$
for $\zz\in U$, $\tau\in\HS$.
\end{thm}

The proof uses the interpretation of indefinite theta series as components 
of triple Fukaya products on a symplectic torus and the
$A_{\infty}$-identity connecting $m_2$ and $m_3$. These products
were first defined by Fukaya in \cite{F}, and in a slightly more
general form by Kontsevich in \cite{K}. 
Basically, the triple products we need correspond to configurations
of four lines with rational slopes on $\R^2$: they are defined
as sums of exponents of areas of the series of quadrangles attached
to such a configuration. These series are always given by some indefinite
theta series as above. It turns out that when two of the four lines
are parallel (i.e. one has a {\it trapezoid} configuration) then
the corresponding quadratic form splits over $\Q$. Now 
using $A_{\infty}$-constraints one can express any
triple Fukaya product on a torus in terms of triple
products corresponding to trapezoid configurations, hence the above
theorem.

The following theorem provides examples of modular indefinite theta series
for which the above principle holds (in other words, the series in 
this theorem correspond to some univalued Massey products on elliptic curve).
Let us introduce a special notation for the summation pattern used for indefinite
theta series: for $S\sub\Q^2$ we denote
$$\sum_{(m,n)\in S}^{indef}a_{m,n}:=\sum_{(m,n)\in S, m\ge 0, n\ge 0}a_{m,n}-
\sum_{(m,n)\in S, m<0, n<0}a_{m,n}.$$

\begin{thm}\label{modser} 
Let $a,b,c,p$ be positive integers such that $a|b$, $c|b$,
$p|(b/a+1)$, $p|(b/c+1)$ and $D=b^2-ac>0$.
Let also $s_1,s_2$ be odd integers and let $r$ be a residue in $\Z/p\Z$.
Then the series
$$q^{\frac{p^2ac(2bs_1s_2-as_1^2-cs_2^2)}{8D}}\cdot
\sum_{(m,n)\in\Z^2, m\equiv n\equiv r(p)}^{indef}
(-1)^{\frac{n-m}{p}}q^{bmn+a\frac{m^2+mps_1}{2}+
c\frac{n^2+nps_2}{2}}$$ 
is a (meromorphic at cusps)
modular form of weight $1$ with respect to some congruenz-subgroup
of $\SL_2(\Z)$. More precisely, it can be written as a ratio
$(\sum_{i=1}^kP_i/Q_i)/(\sum_{j=1}^lR_j/S_j)$, where
$P_i$, $Q_i$, $R_j$ and $S_j$ are polynomials in theta functions
(with characteristics) of degrees
$2(D-2)$, $4(D-2)$, $2(D-1)$ and $4(D-1)$ respectively.
\end{thm}

In the particular case $a=c=p=1$, $r=0$, $b>1$ we get
the series
$$q^{\frac{2bs_1s_2-s_1^2-s_2^2}{8(b^2-1)}}\cdot
\sum_{(m,n)\in\Z^2}^{indef}
(-1)^{m+n}q^{bmn+\frac{m^2+ms_1}{2}+
\frac{n^2+ns_2}{2}}$$
considered in \cite{KP}. These series (for various $s_1,s_2$)
coincide with the string functions of
highest weight modules over $A_1^{(1)}$ of level $b-1$, multiplied
by $\eta^3$, where $\eta$ is the Dedekind eta-function 
(see section \ref{string} for details).
As was shown in \cite{KP} they are equal to Hecke's indefinite theta series
of certain quadratic modules. We generalize this observation in the following 
theorem.

\begin{thm}\label{hecke-thm}
The series considered in Theorem \ref{modser} is equal to
$$N\cdot\Th^H_{\La,Q;\cc}(p^2\tau)$$
for some non-zero integer $N$, where
$$\La=\{(m,n)\in\Z^2: n\equiv(\frac{b}{c}+1)m\mod(2)\},$$
$$\frac{1}{2}Q(m,n)=cn^2-\frac{D}{c}m^2,$$
$$\cc=(\frac{r}{p}+\frac{acs}{2D},\frac{1}{2}).$$
Here $s=\frac{b}{a}s_2-s_1$, so it can be an arbitrary integer
such that $s\equiv \frac{b}{a}+1\mod(2)$. 
\end{thm}

Theorem \ref{modser} is a consequence of Theorem \ref{linsys} below.
We use the following notation:
for a subgroup $I\subset\Z$ and an element $c\in\Q/I$ we denote
$$\theta_{I,c}(z,\tau)=\sum_{m\in c+I}
\exp(\pi i\tau m^2+2\pi i mz).$$
Sometimes we abbreviate $\theta_{I,0}$ to $\theta_I$ and
$\th_{\Z,c}$ to $\th_c$.
 
\begin{thm}\label{linsys} 
Let $d_0$, $d_1$, $d_2$ and $d$ be positive integers
satisfying $d_0+d=d_1+d_2$, $d_1<d$, $d_2<d$.
Let $Q$ be the following quadratic form:
$$Q(m,n)=\frac{1}{d}(d_1(d-d_1)m^2+2d_1d_2mn+d_2(d-d_2)n^2).$$
Fix $\tau$ in the upper half-plane and 
let $(x_1,x_2)$ be a pair of complex numbers satisfying
\begin{equation}\label{disjoint}
x_1+x_2+\frac{1}{2d_1}-\frac{1}{2d_2}\not\in\frac{1}{N}\Z+\Z\tau,
\end{equation}
where $N$ is the least common multiple of $d_1$ and $d_2$.
For a collection of complex numbers $(c_k, k\in\Z/d\Z)$ and an integer $l$,
$0\le l<d_0$, let us consider the series
$$F_l=-\sum_{(m,n)\in\Z^2}^{indef} c_{d_1m-d_2n+l}a_{m,n,l},$$
where 
$$a_{m,n,l}=\exp(\pi i\tau Q(m+\frac{l}{d_0},n-\frac{l}{d_0})+
2\pi i[d_1x_1(m+\frac{l}{d_0})+d_2x_2(n-\frac{l}{d_0})]).$$
Assume that we have
\begin{equation}\label{rowcol}
\sum_{n\in\Z}c_{d_1m_0-d_2n+l}a_{m_0,n,l}=\sum_{m\in\Z}
c_{d_1m-d_2n_0+l}a_{m,n_0,l}=0
\end{equation}
for all $m_0,n_0\in\Z$ and all $l$.
Then $F_l$ are uniquely determined from the linear system of equations
$$\sum_{l\in\Z/d_0\Z}D_{k,l}F_l=c_k$$
for $k\in\Z/d\Z$, where
\begin{align*}
&D_{k,l}=\frac{1}{d_1 i\eta^3(d_1\tau)}\times\\
&\sum_{a\in\Z/d_1\Z}
(-1)^a
\frac{\th_{d_0\Z,l+\frac{d_0}{2}}(\frac{(d_2-d)x_1+d_2x_2}{d_0}+
\frac{2a+1}{2d_1},\frac{\tau}{d_0})
\th_{d\Z,-k+\frac{d}{2}}(x_1+\frac{2a+1}{2d_1},\frac{\tau}{d})}
{\th_{d_2\Z,\frac{d_2}{2}}(x_1+x_2+\frac{2a+1}{2d_1},\frac{\tau}{d_2})},
\end{align*}
$\eta(\tau)=q^{1/24}\cdot\prod_{n\ge 1}(1-q^n)$ is the Dedekind
eta-function.
\end{thm}

The proof of this theorem is also based on the interpretation of indefinite
theta series as matrix coefficients of triple Fukaya products on a torus.
We use the homological mirror symmetry for elliptic curve (see \cite{P-hmc})
to relate these products to Massey products in the derived category of 
an elliptic curve. 

In the case $d_0=1$ Theorem \ref{linsys} gives an explicit formula for
the series $F_0$. In particular, we obtain the following interesting identities 
between $q$-series. 

\begin{corint}\label{identities-cor} One has 
\begin{equation}\label{id1}
q^{\frac{1}{12}}\sum_{(m,n)\in\Z^2}^{indef}
(-1)^{m+n}q^{2mn+\frac{m^2+m}{2}+\frac{n^2+n}{2}}=
\eta(\tau)^2,
\end{equation}
\begin{equation}\label{id2}
\begin{array}{l}
\sum_{(m,n)\in\Z^2,m\equiv n+1(2)}^{indef}
(-1)^{\frac{m+n-1}{2}}q^{\frac{m^2+6mn+3n^2}{2}}=
\frac{\eta^3(2\tau)\theta_{\frac{1}{2}}(\frac{1}{4},3\tau)}
{\theta(\frac{1}{2},4\tau)\theta_{\frac{1}{2}}(\frac{1}{4},\tau)}=\\
q^{\frac{1}{2}}\cdot\prod_{n\ge 1}(1+q^n)(1-q^{2n})(1-q^{3n})(1+q^{6n}),
\end{array}
\end{equation}
\begin{equation}\label{id3}
\begin{array}{l}
\sum_{m\in\Z+\frac{1}{2},n\in\Z}^{indef}
(-1)^{m+n-\frac{1}{2}}q^{\frac{m^2+6mn+3n^2}{2}}=
\frac{\eta^3(2\tau)\theta(\frac{1}{4},3\tau)}
{\theta_{\frac{1}{4}}(0,4\tau)\theta(\frac{1}{4},\tau)}=\\
q^{\frac{1}{8}}\cdot\prod_{n\ge 1}(1-q^n)(1+q^{2n})(1-q^{3n})(1+q^{6n-3}),
\end{array}
\end{equation}
\begin{equation}\label{id4}
\begin{array}{l}
\sum_{m\in\Z,n\in\Z+\frac{1}{2}}^{indef}
(-1)^{m+n-\frac{1}{2}}q^{\frac{m^2+6mn+3n^2}{2}}=
\frac{\eta^3(2\tau)\theta_{\frac{1}{2}}(\frac{1}{4},3\tau)}
{\theta_{\frac{1}{4}}(0,4\tau)\theta_{\frac{1}{2}}(\frac{1}{4},\tau)}=\\
q^{\frac{3}{8}}\cdot
\prod_{n\ge 1}(1-q^n)(1+q^{2n-1})(1-q^{3n})(1+q^{6n}).
\end{array}
\end{equation}
\end{corint}

Identity \eqref{id1} was obtained in \cite{KP} (formula (5.19))
by representation-theoretic means. Three other identities above seem to
be new.

Note that since indefinite theta series are given by alternating sums,
the important problem (raised already by Hecke in \cite{He})
is to determine exactly which of them vanish identically. 
There is a necessary condition (cf. Satz 2 in \cite{He}):
if $\Th^H_{\La,Q;\cc}\neq 0$ then every automorphism of $Q$ preserving
$\La+\cc$ should preserve each component of the cone $Q>0$.
We will prove the following non-vanishing result.

\begin{thm}\label{nonvan}
In the notations of Theorem \ref{modser} consider the series
$$f_{s_1,s_2}=q^{\frac{p^2ac(2bs_1s_2-as_1^2-cs_2^2)}{8D}}\cdot
\sum_{(m,n)\in\Z^2, m\equiv n(p)}^{indef}
(-1)^{\frac{n-m}{p}}
\zeta_{p}^{rm}q^{bmn+a\frac{m^2+mps_1}{2}+
c\frac{n^2+nps_2}{2}},
$$
where $\zeta_p$ is the primitive root of unity of order $p$.
Let us denote by $h$ the greatest common divisor of $b/a+1$, $b/c+1$.
Assume that either
$$\frac{s_1+s_2}{2}\not\in\frac{h}{p}\Z$$
or
$$r\not\in\frac{(h+p)(2b+a+c)}{2hb}+\frac{p}{b}(a\Z+c\Z).$$
Then there exist integers $l_1$ and $l_2$ such that
$$f_{s_1+2\frac{h}{p}l_1,s_2+2\frac{h}{p}l_2}\neq 0$$
\end{thm}

Together with Theorem \ref{hecke-thm} this leads to the following

\begin{corint} In the notations of Theorems \ref{hecke-thm}
and \ref{nonvan} consider the collection of characteristics
$$\cc(t)=(\frac{(b+a)c}{2D}+\frac{act}{D},\frac{1}{2})$$
where $t\in\Z$. Then for every non-zero residue $\ov{t}$ modulo
$\frac{h}{p}$ there exists $t\equiv\ov{t} (\frac{h}{p})$ such that
$\Th^H_{\La,Q;\cc(t)}\neq 0$.
\end{corint}

Using similar techniques we will show in Theorem \ref{det-Jac-thm}
that certain functions of the form
$$\sum_i c_i(\zz)\Theta_{\La,Q,C}(\zz+\vv_i\tau+\ww_i,\tau)$$
are meromorphic Jacobi forms (for a congruenz-subgroup of $\SL_2(\Z)$
and with respect to some quadratic form $Q'$ on a sublattice of $\La$)
in the sense of the definition given by L.~G\"ottsche and D.~Zagier in
\cite{GZ}.
Here $(\vv_i,\ww_i)$ is a collection of vectors in $\La_{\Q}$,
the coefficients $c_i(\zz)$ are products of elliptic functions of some
linear functionals of $\zz$.
In the case of a split form $Q$ we obtain the following result.

\begin{thm}\label{split-Jac-thm} The series
\begin{equation}\label{Jacform-ser}
u^{\frac{s}{2a}}\sum_{n\in\Z}\frac{(-1)^nq^{\frac{n^2+sn}{2}}}{1-q^{an}u},
\end{equation}
where $a$ is a positive integer and $s$ is an odd integer, defines
a meromorphic Jacobi form (here we use multiplicative variables
$q=\exp(2\pi i\tau)$ and $u=\exp(2\pi i z)$).
\end{thm}

In fact, Theorem \ref{linsys} gives an explicit (although complicated) 
rational expression for the series \eqref{Jacform-ser} in terms
of theta functions. In the case $a=1$ such an
expression is well known (see \cite{TM}, Section 486, or \cite{KP} (5.26),
or \cite{P-ap}). In the case $a=2$ this expression is given by formula
\eqref{a2-eq}.

Considering higher Fukaya products $m_k$ with $k\ge 4$ on a symplectic
torus one
still gets some indefinite theta series corresponding to lattices with
quadratic forms of signature $(1,k-2)$. However, quadratic forms associated
with configurations of $k+1$ lines depend on $k$ parameters while
a general quadratic form on a lattice of rank $k-1$ has
$(k-1)k/2$ coefficients. Hence, for $k\ge 4$ not every indefinite
series associated with a quadratic form of signature $(1,k-2)$ 
comes from a Fukaya product.

Here is the plan of the paper. In section 1 we explain the relation
between indefinite theta series of signature $(1,1)$ and the Appell's
function
(\ref{kappa}). Section 2 contains the definition of Fukaya products
on a symplectic torus and the computation of double and
triple Fukaya products in terms of the usual theta functions and
indefinite theta series respectively. In section 3 we prove an auxiliary
surjectivity result about the products $m_2$ in the Fukaya category of
a symplectic torus. In section 4 we prove Theorem \ref{main}
and in section 5 we illustrate it by an explicit example. Section 6
is devoted to the definition and computation of Massey products of morphisms
between line bundles on elliptic curve. In section 7 we give examples
of modular indefinite theta series, proving in particular Theorems
\ref{modser}, \ref{hecke-thm}, \ref{linsys}, \ref{nonvan}, \ref{split-Jac-thm}
and Corollary \ref{identities-cor}.

\noindent
{\it Acknowledgment}. During the preparation of this paper
I benefited from conversations with B. Gross, M. Kontsevich and 
D. Zagier. I am grateful to V. Kac for the reference to Hecke's works and
to \cite{KP}. Part of this paper was written during my visit to
Max-Planck-Institut f\"ur Mathematik. I'd like to thank the Institute
for its hospitality. This work was partially supported by the NSF grant.

\section{Indefinite theta series of signature $(1,1)$}\label{thetasec}

Let $\La$ be a rank 2 lattice equipped with a $\Q$-valued quadratic form $Q$,
$\tau$ be an element in the upper-half plane. We assume that $Q$ has
signature $(1,1)$ and fix a rational open cone $C\in\La_{\R}$ such that
$Q|_C>0$. We will use the following notation for indefinite theta series with
characteristics associated with $(\La,Q,C)$: for an element 
$\cc\in\La_{\Q}/\La$ we set
\begin{equation}\label{char}
\Theta_{\La,Q,C;\cc}(\zz,\tau)= 
\sum_{\nn\in \cc+\La:\nn+\a(\zz)\in C}\sign(\nn+\a(\zz))
\exp(\pi i\tau Q(\nn)+2\pi i \nn\cdot \zz)
\end{equation}
In other words, we have
$$\Theta_{\La,Q,C;\cc}(\zz,\tau)=\exp(\pi i\tau Q(\cc)+2\pi i \cc\cdot \zz)
\Theta_{\La,Q,C}(\zz+\tau \cc,\tau).$$

The following identities follow immediately from the definition:
$$\Th_{N\La,Q,C,\cc}(\zz,\tau)=\Th_{\La,Q,C,\frac{\cc}{N}}(N\zz,N^2\tau),$$
$$\Th_{\La,NQ,C,\cc}(\zz,\tau)=\Th_{\La,Q,C,\cc}(N\zz,N\tau)$$
where $N>0$ in an integer,
$$\Th_{\La,Q,C,\cc}(\zz,\tau)=\sum_{\nn\in\La/\La'}
\Th_{\La',Q,C;\cc+\nn}(z,\tau)$$
for any sublattice $\La'\sub\La$.
Since $N\La\sub\La'$ for some $N$ we can also use these formulas 
to express $\Th_{\La',Q,C}$ in terms of $\Th_{\La,Q,C}$.  

On the other hand, since the cone $C$ is rational we can choose coordinates
in such a way that $C=\{(x_1,x_2)\in\R^2: x_1x_2>0\}$, $x_i>0$ in $C^+$
and $\La$ is a lattice in $\R^2$ commensurable with $\Z^2$.
Now the condition $Q|_C>0$ and the requirement that the signature of $Q$ is
$(1,1)$ mean that $Q(x_1,x_2)=a_{11}x_1^2+2a_{12}x_1x_2+a_{22}x_2^2$
where $a_{ii}\ge 0$ and $D=a_{12}^2-a_{11}a_{22}>0$.

Now let us consider the case when $Q$ splits into a product of linear
forms over $\Q$. Then by additivity of $\Th$ in $C$
it suffices to consider the case when $Q$ vanishes on one of the lines
forming the boundary of $C$.
Then we can choose coordinates in such a way that
$C=\{(x_1,x_2)\in\R^2: x_1x_2>0\}$, $x_i>0$ in $C^+$,
$Q(n_1,n_2)=an_1(n_1+2n_2)$ for some $a\in\Q$, $\La$ is a lattice
commensurable with $\Z^2$. Rescaling $\tau$ and $\zz$ (rationally)
we can assume
that $a=1$. Also it suffices to consider the case $\La=\Z^2$.
Then for any $\cc=(c_1,c_2)$ we have
\begin{align*}
&\Th_{\Z^2,Q,C,\cc}(\zz,\tau)=
\sum_{\nn\in\Z^2+\cc,(n_1+\a(z_1))(n_2+\a(z_2))>0}\sign(n_1+\a(z_1))\times\\
&\exp(\pi i\tau n_1(n_1+2n_2)+2\pi i(n_1z_2+n_2z_1)+2\pi i n_1z_1).
\end{align*}
We can split this sum in two pieces and sum the geometric
progression in $n_2$ in each of them:
\begin{align*}
&\Th_{\Z^2,Q,C,\cc}(\zz,\tau)=\\
&\sum_{n_1\in\Z+c_1,n_1+\a(z_1)>0}\exp(\pi i\tau n_1^2+2\pi i n_1(z_1+z_2))
\sum_{n_2\in\Z_{\ge 0}+n_2^0}\exp(2\pi i n_2(\tau n_1+z_1))-\\
&\sum_{n_1\in\Z+c_1,n_1+\a(z_1)<0}\exp(\pi i\tau n_1^2+2\pi i n_1(z_1+z_2))
\sum_{n_2\in\Z_{\le 0}+n_2^0-1}\exp(2\pi i n_2(\tau n_1+z_1))=\\
&\sum_{n_1\in\Z+c_1,n_1+\a(z_1)>0}
\frac{\exp(\pi i\tau n_1^2+2\pi i n_1(z_1+z_2)+2\pi i n_2^0(\tau n_1+z_1))}
{1-\exp(2\pi i(\tau n_1+z_1))}-\\
&\sum_{n_1\in\Z+c_1,n_1+\a(z_1)<0}
\frac{\exp(\pi i\tau n_1^2+2\pi i n_1(z_1+z_2)+2\pi i (n_2^0-1)(\tau n_1+z_1))}
{1-\exp(-2\pi i(\tau n_1+z_1))}=\\
&\sum_{n_1\in\Z+c_1}
\frac{\exp(\pi i\tau n_1^2+2\pi i n_1(z_1+z_2)+2\pi i n_2^0(\tau n_1+z_1))}
{1-\exp(2\pi i(\tau n_1+z_1))}
\end{align*}
where $n_2^0$ is the minimal $n_2\in\Z+c_2$ such that $n_2+\a(z_2)>0$.
Hence, we derive the following formula:
\begin{equation}\label{appell}
\Th_{\Z^2,Q,C,\cc}(\zz,\tau)=
\exp(2\pi i n_2^0z_1)\kappa_{c_1}(z_1,(1-n_2^0)\tau-z_1-z_2;\tau).
\end{equation}

\section{Fukaya category of a torus}

\subsection{Definition}\label{Fukdef-sec}

Let us recall the definition of the Fukaya $A_{\infty}$-category of the 
torus $\R^2/\Z^2$ with the (complexified) symplectic form 
$-2\pi i\tau dx\wedge dy$ where $\tau$ is an element of the upper half-plane
(for more details see \cite{P-hmc}).
More precisely, it is not quite an $A_{\infty}$-category since morphisms are only
defined for transversal configurations of objects, however, the axiomatics
can be changed appropriately (see \cite{KS}, sec. 4.3)
Also, we will need only the subcategory $\FF_s$ which is described as
follows. The objects of $\FF_s$ are pairs $(L,t)$
where  $L\subset\R^2$ is a non-vertical 
line with rational slope considered modulo translations
by $\Z^2$, $t$ is a real number.
Morphisms between two such objects $(L_1,t_1)$ and $(L_2,t_2)$ are 
defined only if $L_1\neq L_2\mod \Z^2$. In this case
$\Hom((L_1,t_1),(L_2,t_2))=\Hom(L_1,L_2)$ is a $\C$-vector space 
with the basis $[P]$ enumerated
by points $P\in (L_1+\Z^2)\cap (L_2+\Z^2)$ modulo $\Z^2$ 
(the numbers $t_i$ will
play a role only in the definition of compositions).
Let $\la_i$ be the slope of the line $L_i$ ($i=1,2$). Then 
$\Hom(L_1,L_2)\neq 0$ only if $\la_1\neq\la_2$.  
This space has grading $0$ if $\la_1<\la_2$ and grading $1$ if 
$\la_1>\la_2$.
By definition the differential $m_1$ is zero. The compositions $m_k$ 
for $k\ge 2$ are (partially) defined as follows.
Let $L_0,L_1,\ldots,L_k$ be the set of lines in $\R^2$ with slopes 
$\la_0,\la_1,\ldots,\la_k$. 
Assume that
the images of $L_i$ in $\R^2/\Z^2$ form a {\it transversal configuration},
i.e., no three of them intersect in one point. 
For every $i=0,\ldots,k-1$ let
$d_i$ be the grading of $\Hom(L_i,L_{i+1})$. The composition
$$m_k:\Hom(L_0,L_1)\otimes\ldots\otimes
\Hom(L_{k-1},L_k)\ra\Hom(L_0,L_k)$$
is non-zero only if $\sum_{i=0}^{k-1}d_i-k+2$ is equal to the degree
of $\Hom(L_0,L_k)$. 
Let $P_{i,i+1}$ be some intersection points of $L_i$
and $L_{i+1}$ modulo $\Z^2$. Then
\begin{align*}
&m_k([P_{0,1}],[P_{1,2}],\ldots,[P_{k-1,k}])=\\
&\sum_{P_{0,k},\Delta}\pm
\exp\left(2\pi i\tau\cdot\int_{\Delta}dx\wedge dy+
2\pi i\sum_{j\in\Z/(k+1)\Z}(x(p_j)-x(p_{j-1}))t_j\right)[P_{0,k}]
\end{align*}
where the sum is taken over points of intersections $P_{0,k}$ of
$L_0$ with $L_k$ modulo $\Z^2$ and over all $(k+1)$-gons $\Delta$ 
(considered up to traslation by $\Z^2$) with vertices
$p_i\equiv P_{i,i+1} \mod\Z^2$, $i\in\Z/(k+1)\Z$,
such that the edge $[p_{i-1},p_i]$ belongs to $L_i+\Z^2$. We also require
that the path formed by the edges $[p_0,p_1], [p_1,p_2], \ldots, [p_k,p_0]$
goes in the clockwise direction. 
The sign in the RHS is ``plus" if $k$ is even and is equal to
the sign of $x(p_0)-x(p_k)$ if $k$ is odd.

The $A_{\infty}$-constraint
we are going to use is
\begin{eqnarray}\label{ainfty}
m_3(m_2(a_1,a_2),a_3,a_4)-m_3(a_1,m_2(a_2,a_3),a_4)+m_3(a_1,a_2,m_2(a_3,a_4))
=\nonumber \\
= m_2(m_3(a_1,a_2,a_3),a_4)+(-1)^{\deg(a_1)}m_2(a_1,m_3(a_2,a_3,a_4)),
\end{eqnarray}
where $a_1,\ldots,a_4$ are composable morphisms between $5$ objects in $\FF_s$
forming a transversal configuration.
Below we will often abbreviate $m_2(a,b)$ to $ab$.  

\subsection{Double products and vector bundles on elliptic curves}
\label{double}

Since $m_1=0$ the composition $m_2$ is associative, so we can consider
the category $\FF_s$ with $m_2$ as a usual category.
It was shown in \cite{PZ} that the obtained category is equivalent to the
category of stable vector bundles on the elliptic curve
$E=\C/\Z+\Z\tau$, where morphisms between vector bundles $V_1$ and $V_2$
are elements of the graded vector space $\oplus_i \Ext^i(V_1,V_2)$
(in fact, we showed in \cite{PZ} how to extend this equivalence to the whole 
derived category of coherent sheaves on elliptic curve, but we don't need
this extension here). The construction of this equivalence (which we
recall below) is based on the observation due to M.~Kontsevich 
that the Fukaya product $m_2$ on a torus is given essentially
by theta functions. Here is a more precise statement.

For a pair $(\la,y)$ where $\la\in\Q$, $y\in\R$, let us denote
by $L(\la,y)$ the line in $\R^2$ given by
$$L(\la,y)=\{(t,\la t-y), t\in\R\}.$$
Now let $L_i=L(\la_i,y_i)$, $i=0,1,2$,
be lines in $\R^2$ with distinct slopes $\la_i$. 
Let us denote
\begin{eqnarray}\label{yij}
y_{ij}=\frac{y_j-y_i}{\la_j-\la_i}, \nonumber\\
y'_{ij}=\frac{\la_iy_j-\la_jy_i}{\la_j-\la_i}\nonumber. 
\end{eqnarray}
The lines $L_i$ and $L_j$ intersect at the point
$$P_{ij}(y_i,y_j)=(y_{ij},y'_{ij}).$$
Note that
if we change $y_j$ by $y_j+m\la_j+n$ where $m,n\in\Z$,
the new line $L(\la_j,y_j+m\la_j+n)$ is a shift of $L_j$ by an integer
vector. Thus, the new point of intersection $P_{ij}(y_i,y_j+m\la_j+n)$
still belongs to $L_i\cap(L_j+\Z^2)$. One has 
$P_{ij}(y_i,y_j+m\la_j+n)\equiv P_{ij}(y_i,y_j)\mod\Z^2$ if and
only if $(m,n)\in \La(\la_i,\la_j)$ where
$$\La(\la_i,\la_j)=\{(m,n)\in\Z^2:\ 
\frac{m\la_j+n}{\la_j-\la_i}\in I_{\la_i}\}$$ 
where for every $\la\in\Q$ we denote 
\begin{equation}\label{Ila}
I_{\la}=\{n\in\Z: \la n\in\Z\}.
\end{equation} 
Thus, we have the following basis in $\Hom(L_i,L_j)$:
$$[P_{ij}(y_i,y_j+m\la_j+n)], (m,n)\in\Z^2/\La(\la_i,\la_j).$$
Instead of shifting $y_j$ we could also shift $y_i$ and get a different
indexing of intersection points modulo $\Z^2$. However, this indexing 
is related to the previous one by the formula
$$P_{ij}(y_i-m\la_i-n,y_j)=P_{ij}(y_i,y_j+m\la_j+n).$$
Note also that we have $\La(\la_i,\la_j)=\La(\la_j,\la_i)$, so
changing the order of lines we would get essentially the same indexing.

Assume that 
$\deg\Hom(L_0,L_1)+\deg\Hom(L_1,L_2)=\deg\Hom(L_0,L_2)$.
Let $t_i$, $i=0,1,2$, be some real numbers. Then we can consider objects
$(L_i,t_i)$ in Fukaya category.
An easy computation shows that 
\begin{eqnarray}\label{m2P}
m_2([P_{01}(y_0,y_1)],[P_{12}(y_1,y_2)])=\qquad\qquad\qquad\nonumber\\
\sum_{n\in I_{\la_1}} \exp(\pi i \tau p(v_1+n)^2-
2\pi i p(v_1+n)w_1)[P_{02}(y_0,y_2+n\la_2-n\la_1)]
\end{eqnarray}
where 
$$p=p(\la_0,\la_1,\la_2)=\frac{(\la_2-\la_1)(\la_1-\la_0)}{(\la_2-\la_0)},$$
$v_1=y_{12}-y_{01}$, $w_1=t_{12}-t_{01}$, 
$$t_{ij}=\frac{t_j-t_i}{\la_j-\la_i},$$
$I_{\la_1}$ is defined by (\ref{Ila}).
Note that the class of the point $[P_{02}(y_0,y_2+n\la_2-n\la_1)]$ modulo $\Z^2$ 
depends only on the class
of $n$ modulo the following subgroup
$$I_{\la_0,\la_1,\la_2}=I_{\la_1}\cap\frac{\la_2-\la_0}{\la_2-\la_1}I_{\la_0}.
$$
The matrix coefficients of the above product are given by values of 
elliptic functions
at $(\tau v_1-w_1,\tau)$ times the non-holomorphic factor
$\exp(\pi i\tau p v_1^2-2\pi i p v_1w_1)$. 
One can get rid of this factor by rescaling the bases in $\Hom(L_i,L_j)$
appropriately. Namely, we set
\begin{equation}\label{eij}
e_{ij}(m,n)=e_{y_i,y_j}(m,n)=
\exp(\pi i\tau(\la_i-\la_j)y_{ij}^2-2\pi i(t_i-t_j)y_{ij})
[P_{ij}(y_i,y_j+m\la_j+n)]
\end{equation}
where $(m,n)\in\Z^2/\La(\la_i,\la_j)$.
Then the above formula is equivalent to 
$$m_2(e_{01}(0,0),e_{12}(0,0))=
\sum_{n\in I_{\la_1}/I_{\la_0,\la_1,\la_2}} 
\theta_{I_{\la_0,\la_1,\la_2},n}(p(v_1\tau-w_1),p\tau)
e_{02}(n,-n\la_1)
$$
where we use the notation $\theta_{I,c}$ from the introduction.
Changing $y_i$'s appropriately in the formula (\ref{m2P}) we derive
a more general formula
\begin{eqnarray}\label{m2e}
m_2(e_{01}(a,b),e_{12}(c,d))=\qquad\qquad\qquad\qquad\nonumber\\
\sum_{n\in I_{\la_1}/I_{\la_0,\la_1,\la_2}} 
\theta_{I_{\la_0,\la_1,\la_2},u+n}(p(v_1\tau-w_1),p\tau)
e_{02}(a+c+n,b+d-n\la_1)
\end{eqnarray}
where 
$$u=\frac{c\la_2+d}{\la_2-\la_1}-\frac{a\la_0+b}{\la_1-\la_0}.$$
The corresponding coefficients will be holomorphic in 
$v_1\tau-w_1$. 
The associativity condition for $m_2$ is equivalent to the classical
addition formulas for theta-functions.

The equivalence with the category of stable bundles on elliptic curve
$E_{\tau}=\C/\Z+\tau\Z$ is constructed in \cite{PZ} as follows.
First let us consider the subcategory in $\FF_s$ formed
by lines with integer slopes. To an object of this subcategory
$(L(\la,y),t)$ where $\la\in\Z$, $y,t\in\R$, we associate
the line bundle $t^*_{y\tau-t}\LL\otimes\LL^{\otimes(\la-1)}$
where $\LL=\LL_{\tau}$ is the line bundle of degree $1$ on $E_{\tau}$
such that $\theta(z,\tau)$ is a holomorphic section of $\LL_{\tau}$,
$t_z:E_{\tau}\ra E_{\tau}$ denotes the translation by $z$.
Assume that we have 
two such objects $(L_i,t_i)$, $i=1,2$, where $L_i=L(\la_i,y_i)$,
$\la_i\in\Z$, $\la_1<\la_2$. Then 
we identify $\Hom(L_1,L_2)$ with the space of morphisms between the
corresponding line bundles by sending the basis elements
$e_{12}(0,k)$, $k\in\Z/(\la_2-\la_1)\Z$, defined by (\ref{eij}) to
the functions
$$\theta_{(\la_2-\la_1)\Z,k}(z+y_{12}\tau-t_{12},\frac{\tau}{\la_2-\la_1})$$
regarded as holomorphic sections of the line bundle
$$(t^*_{y_1\tau-t_1}\LL\otimes\LL^{\otimes(\la_1-1)})^*\otimes
(t^*_{y_2\tau-t_2}\LL\otimes\LL^{\otimes(\la_2-1)})\simeq
t^*_{y_{12}\tau-t_{12}}\LL^{\otimes(\la_2-\la_1)}.$$
The fact that this map respects $m_2$ follows from addition formulas
for theta functions. To extend this equivalence to all lines and
all stable bundles we use isogenies.
For every positive integer $r$ consider the natural isogeny of degree $r$
$$\pi_r:E_{r\tau}\ra E_{\tau}.$$ 
Then we have the natural functors $\pi_{r*}$ and $\pi_r^*$ 
between the categories of bundles on $E_{\tau}$ and $E_{r\tau}$.
We complete the construction of our equivalence by requiring that these 
functors correspond
to the obvious functors $\pi_{r*}$ and $\pi_r^*$ between the corresponding
Fukaya categories (see \cite{PZ} for details). One also has to identify
morphisms of degree $1$ in both categories. 
For this one has to fix a non-zero holomorphic
$1$-form on $E_{\tau}$ and use the isomorphisms $\Hom(V_1,V_2)^*\simeq
\Ext^1(V_2,V_1)$ (where $V_1$ and $V_2$ are vector bundles on $E_{\tau}$)
induced by Serre duality together with the obvious isomorphisms
$\Hom^0(L_1,L_2)^*\simeq\Hom^1(L_2,L_1)$ in the Fukaya category.

\subsection{Triple products and indefinite theta series}
\label{triple}

Consider $4$ lines 
$(L_i=L(\la_i,y_i), i\in\Z/4\Z)$ where $\la_i\in\Q$,
$y_i\in\R$. As before, 
we assume that the corresponding circles in $\R^2/\Z^2$ form
a transversal configuration, in particular, 
the lines $L_i$ are distinct modulo $\Z^2$ and 
$\la_i\neq\la_{i+1}$ for $i\in\Z/4\Z$.
Let $t_i, i\in\Z/4\Z$ be some real numbers, then $(L_i,t_i)$ are objects of the
Fukaya category.  We are going to compute 
the Fukaya triple product $m_3([P_{01}],[P_{12}],[P_{23}])$,
where $P_{i,i+1}:=P_{i,i+1}(y_i,y_{i+1})$ for $i=0,1,2$.
This product is zero unless the following equality is
satisfied:
\begin{equation}\label{degree}
\sum_{i=0}^2\deg\Hom(L_i, L_{i+1})=\deg\Hom(L_0,L_3)+1.
\end{equation}
Following the definition we have to consider all
quadrangles (up to translation by $\Z^2$)
with vertices $p_i$
such that for every $i$ the vector $p_i-p_{i-1}$ has slope $\la_i$, 
$p_i\equiv P_{i,i+1}\mod \Z^2$, for $i=0,1,2$, and
the piecewise linear path $[p_0,p_1,p_2,p_3]$ goes in the clockwise direction.
First of all, it is easy to check that the condition (\ref{degree}) implies
that all such quadrangles are convex. Secondly, the condition on the
orientation of the path is equivalent to the system of inequalities  
\begin{equation}\label{ineq-det}
\det(p_{i+1}-p_i,p_i-p_{i-1})>0
\end{equation}
where $i\in\Z/4\Z$.
These quadrangles (considered up to translations by $\Z^2$)
can be parametrized by elements of a rank-$2$
lattice. Namely, let
$$\La=\La(\la_0,\ldots,\la_3)=\{\nn=(n_0,\ldots,n_3)\in\Q^4:\ 
\sum n_i=\sum\la_in_i=0, n_1\in I_{\la_1}, n_2\in I_{\la_2}\}.$$
Then writing 
$$p_i-p_{i-1}=x_i(1,\la_i)$$
for $i\in\Z/4\Z$ we obtain the vector
$\xx=(x_0,\ldots,x_3)$ in $\La_{\R}$. 
The inequalities (\ref{ineq-det}) become
\begin{equation}\label{ineq-lat}
(\la_i-\la_{i+1})x_ix_{i+1}>0
\end{equation}
for $i\in\Z/4\Z$. On the other hand, setting
$$P_{i,i+1}-P_{i-1,i}=v_i(1,\la_i)$$
we obtain the vector $\vv=(v_0,\ldots,v_i)\in\La_{\R}$
(note that $v_i=y_{i,i+1}-y_{i-1,i}$).
Now the conditions $p_i\equiv P_{i,i+1}\mod\Z^2$ for $i=0,1,2$ imply that
$\xx-\vv$ belongs to $\La$. Conversely, given an element
$\nn=(n_0,\ldots,n_3)\in\La$ we have the corresponding quadrangle
$\Delta(\nn)$ with vertices $p_i$ such that $p_0=P_{0,1}$ and
$p_i-p_{i-1}=(v_i+n_i)(1,\la_i)$. Fixing the fourth
vertex $p_3$ modulo $\Z^2$ is equivalent to choosing $\nn$ in a fixed
coset modulo the sublattice $\La^+\subset\La$ defined as follows:
$$\La^+=\La^+(\la_0,\ldots,\la_3)=\{\nn=(n_0,\ldots,n_3)\in\Z^4:\ 
\sum n_i=\sum\la_in_i=0, \la_in_i\in\Z\}.$$
Thus, the sums in the definition
of the Fukaya coefficients are taken over all elements $\nn$
of a coset of $\La^+$ in $\La$,
such that $\vv+\nn\in C$, where $C\subset\La_{\R}$ is an open subset defined
by inequalities (\ref{ineq-lat}). It is easy to see that
the condition (\ref{degree}) implies that $C$ is a non-empty open cone.
The relation between the vertex $p_3$ and the coset $\nn\in\La/\La^+$
can be found explicitly as follows.
We know that $p_0-p_3=(v_0+n_0)(1,\la_0)$,
and that $p_0\equiv P_{01}$. Hence, $p_3\equiv P_{03}(y_0,y_3+a\la_3+b)$
where $a$ and $b$ are integers satisfying
$$\frac{a\la_3+b}{\la_3-\la_0}\equiv -n_0\mod(I_{\la_0}).$$
It follows that
$$(a,b)\equiv (n_1+n_2,-\la_1n_1-\la_2n_2)\mod\La(\la_0,\la_3).$$
The area of $\Delta(\nn)$ is given by
$$\int_{\Delta(\nn)}dx\wedge dy=
\frac{1}{2}\left(\det(p_1-p_0,p_0-p_3)+\det(p_3-p_2,p_2-p_1)\right)=
\frac{1}{2}Q(\vv+\nn)$$
where $Q$ is the quadratic form on $\La_{\R}$ given
by 
$$Q(\xx)=(\la_0-\la_1)x_0x_1+(\la_2-\la_3)x_2x_3.$$

Finally, we have
$$\sum_{i\in\Z/4\Z}(x(p_i)-x(p_{i-1}))t_i=\sum_i t_ix_i=
-\ww\cdot \xx$$
where $\ww=(w_0,\ldots,w_3)\in\La_{\R}$, $w_i=t_{i,i+1}-t_{i-1,i}$,
$t_{i,j}=\frac{t_j-t_i}{\la_j-\la_i}$, $\xx\cdot \xx'$ is the symmetric
pairing induced by $Q$ (so that $Q(\xx)=\xx\cdot \xx$).

Thus, we have
\begin{align*}
&m_3([P_{01}(y_0,y_1)],[P_{12}(y_1,y_2)],[P_{23}(y_2,y_3)])=\\
&\sum_{\nn\in\La, \vv+\nn\in C}\pm \exp(\pi i\tau Q(\vv+\nn)-
2\pi i \ww\cdot (\vv+\nn))
[P_{03}(y_0,y_3+(n_1+n_2)\la_3-\la_1n_1-\la_2n_2)]
\end{align*}
where the sign is equal to the sign of $v_0+n_0$.
Let us choose $C^{+}$ to be the component of $C$ where $x_0>0$.
Then using the bases $e_{ij}(m,n)$ defined by (\ref{eij}) we can
rewrite the above formula as follows:
\begin{align*}
&m_3(e_{01}(0,0),e_{12}(0,0),e_{23}(0,0))=\\
&\sum_{\nn\in\La/\La^+}\Theta_{\La^+,Q,C;\nn}(\tau \vv-\ww,\tau) 
e_{03}(n_1+n_2,-\la_1n_1-\la_2n_2)
\end{align*}
where $\Theta_{\La^+,Q,C;\nn}$ is the indefinite theta series with 
characteristic $\nn$ defined by (\ref{char}).
Similarly we can compute products of all basis elements:
\begin{eqnarray}\label{m3}
m_3(e_{01}(a,b),e_{12}(c,d),e_{23}(f,g))= \qquad\qquad\qquad\qquad \nonumber\\
\sum_{\nn\in\La/\La^+}\Theta_{\La^+,Q,C;\uu+\nn}(\tau \vv-\ww,\tau) 
e_{03}(a+c+f+n_1+n_2,b+d+g-\la_1n_1-\la_2n_2)
\end{eqnarray}
where
$\uu=(u_0,u_1,u_2,u_3)\in\La_{\Q}$ is the following vector:
\begin{align*}
&\uu=(\frac{a\la_1+b}{\la_1-\la_0}-\frac{(a+c+f)\la_3+b+d+g}{\la_3-\la_0},
\frac{c\la_2+d}{\la_2-\la_1}-\frac{a\la_0+b}{\la_1-\la_0},
\frac{f\la_3+g}{\la_3-\la_2}-\frac{c\la_1+d}{\la_2-\la_1},\\
&\frac{(a+c)\la_0+f\la_3+b+d+g}{\la_3-\la_0}-\frac{f\la_3+g}{\la_3-\la_2}).
\end{align*} 

The $A_{\infty}$-axiom (\ref{ainfty}) can be converted 
into a certain identity for indefinite theta series (and the usual theta
functions which appear from $m_2$).
The explicit formula can be found in the last section of \cite{P}.

It is convenient for explicit computations 
to choose a pair of components $(x_i,x_j)$ as coordinates on $\La_{\R}$ in 
such a way that the cone $C$ defined by inequalities (\ref{ineq-lat})
coincides with the cone $x_ix_j>0$. We will do this in two
particular cases. 

\noindent
1. First assume that $\la_1<\la_0\le\la_2<\la_3$. 
Choose $(x_0,x_1)$ as coordinates in $\La_{\R}$. Then
we have $C=\{\xx:\ x_0x_1>0\}$. The quadratic form in these
coordinates can be written as 
$$Q(\xx)=ax_0^2+2bx_0x_1+cx_1^2,$$
where
$$a=\frac{(\la_2-\la_0)(\la_3-\la_0)}{\la_3-\la_2},$$
$$b=\frac{(\la_2-\la_1)(\la_3-\la_0)}{\la_3-\la_2},$$
$$c=\frac{(\la_2-\la_1)(\la_3-\la_1)}{\la_3-\la_2}.$$
If $\la_0=\la_2$ (in this case we say that this
is a {\it trapezoid} triple product) then $a=0$,
so the form $Q$ splits over $\Q$. 
On the other hand, if we set $\la_0=0$ then 
the corresponding transformation $(\la_1,\la_2,\la_3)\mapsto
(a,b,c)$ is birational. More precisely,
the inverse transformation is given by
$$\la_1=-\frac{D}{b},$$
$$\la_2=\frac{aD}{b(b-a)},$$
$$\la_3=\frac{D}{c-b}$$
where $D=b^2-ac$. The inequalities satisfied by $\la_i$ are equivalent to the
following inequalities for $a,b,c$:
$$-D<0\le a<b<c.$$

\noindent
2. Now assume that $\la_2<\la_0<\la_3<\la_1$. Choose
$(x_0,x_3)$ as coordinates in $\La_{\R}$. Then $C=\{\xx:\ x_0x_3>0\}$
and 
$$Q(\xx)=ax_0^2+2bx_0x_3+cx_3^2,$$
where
$$a=\frac{(\la_1-\la_0)(\la_0-\la_2)}{\la_1-\la_2},$$
$$b=\frac{(\la_1-\la_0)(\la_3-\la_2)}{\la_1-\la_2},$$
$$c=\frac{(\la_1-\la_3)(\la_3-\la_2)}{\la_1-\la_2}.$$
Setting $\la_0=0$ we get a birational transformation coinciding with
the previous one up to permutation of variables and signs. 
So the inverse map is given by
$\la_1=D/(b-c)$, $\la_2=aD/b(a-b)$, $\la_3=D/b$.
The inequalities for $\la_i$ are equivalent to the following inequalities:
$$-D<0<a,c<b.$$

It follows that every indefinite theta series associated with 
rational quadratic form of signature $(1,1)$ appears as a coefficient
of certain Fukaya triple product.
Indeed, let us consider the $\Q$-valued quadratic form
$Q(x,y)=ax^2+2bxy+cy^2$ on $\Z^2$ such that $D=b^2-ac>0$ and $Q$ 
is positive on the cone $C_0=\{(x,y)\in\R^2:\ xy>0\}$.
Assume first that $ac\neq 0$ and $b\neq c$. Then we necessarily have that
$a$, $b$ and $c$ are positive and either $a<b$ or $c<b$. Permuting
the coordinates if necessary we can assume that $a<b$. Then 
either $c>b$ and we are in the
situation of the case 1 above or $c<b$ so we can apply the case 2.
Note that the lattice coming from the configuration
of lines will be commensurable with $\Z^2$, so we can use formulas
of section \ref{thetasec} to relate the indefinite theta series associated with
$(\Z^2,Q,C_0)$ to the corresponding Fukaya triple product.
Furthermore, by rescaling the coordinates
$x$, $y$ we can always achieve that we are in the situation of a given
case. For example, for the proof of Theorem \ref{main} we will use
a rescaling which leads to the case 1.
If $a=0$ and $b\neq c$ we can still apply the formulas of either case 1 or
case 2. If $b=c$ then one of the slopes will be infinite.
The only reason why we didn't include vertical lines in our category
was because they correspond to torsion sheaves on elliptic curves while
we want to deal only with bundles. However, the Fukaya compositions
are well-defined for all lines including vertical, so a slight 
modification of the above computation will work in this case.
Alternatively, we can always rescale the coordinates in such a way that
$b\neq c$ and then apply the above formulas. On the other hand,
we notice that replacing the form $Q$ by $NQ$ for some $N>0$
we can always achive that all the slopes $\la_i$ are integers.

\begin{rem} For computations with the form $Q$ on $\La_{\R}$ defined above
the following formula is useful:
$$\xx\cdot \yy=
(\la_i-\la_{i+1})x_iy_{i+1}+(\la_{i+2}-\la_{i+3})y_{i+2}x_{i+3},$$
for any $i\in\Z/4\Z$, $\xx,\yy\in\La_{\R}$.
\end{rem}

\section{Morphisms of vector bundles on elliptic curves}

We identify an elliptic curve $E$ with its dual by associating
to a point $x\in E$ the line bundle $\PP_x=\O_E(x-e)$ of degree zero on
$E$, where $e\in E$ is the neutral element of the group law. 
For every integer $d$ we denote by $E_d$ the kernel of the 
homomorphism $[d]:E\ra E:x\mapsto dx$.

\begin{prop}\label{surj-pr}
Let $V_1$, $V_2$ and $V_3$ be stable vector bundles on an
elliptic curve $E$. Assume that $V_i$ has rank $r_i$ and degree $d_i$
and that the slopes $\mu_i=d_i/r_i$ satisfy $\mu_1<\mu_2<\mu_3$.
Then there exists an integer $d$ depending only on $(d_i)$ and $(r_i)$
such that the following natural map is surjective
$$\oplus_{x\in E_d}\Hom(V_1,V_2\otimes\PP_x)\otimes\Hom(V_2\otimes\PP_x,V_3)
\ra\Hom(V_1,V_3).$$
\end{prop}

\Pf . Consider the action of $E\times E$ on the category of
vector bundles on $E$, such that a point $(x,y)$ acts as the functor
$T_{(x,y)}:F\mapsto t_x^*F\otimes\PP_y$. 
Then the statement of the theorem can be reformulated as follows:
there exists a finite subgroup $S\subset E\times E$ such that the map
\begin{equation}\label{hommap}
\oplus_{s\in S}\Hom(V_1,T_s(V_2))\otimes
\Hom(T_s(V_2),V_3)\ra\Hom(V_1,V_3)
\end{equation}
is surjective. Indeed, this follows from the fact that
for any $x\in E$ one has
$$t^*_{r_2x}V_2\simeq V_2\otimes\PP_{-d_2x}.$$

Now we claim that in proving the surjectivity of (\ref{hommap}) we can
replace the bundles $V_i$ by $V_i\otimes L$, where $L$ is a line bundle,
or by $\SS(V_i)$ provided that $d_1>0$, where $\SS$ is the Fourier-Mukai
transform (see \cite{M}). Indeed, in the first case this is clear.
In the case of the Fourier-Mukai transform this follows from
the fact that $\SS$ interchanges translations with tensoring by line
bundles of degree zero.

Using these two operations (tensoring with a line bundle and the Fourier-Mukai
transform) we can make $V_1=\O_E$. Next we want to reduce the proof to the
case when $V_2$ is a line bundle. Indeed, assume that in this case the
assertion is true. Then consider an isogeny $\pi:E'\ra E$ of degree $r_2$
and a line bundle $L$ on $E'$ such that
$\pi_*L\simeq V_2$ (such $\pi$ and $L$ always exist).
By assumption the statement is true for the triple 
$(\O_{E'},L,\pi^*V_3)$ on $E'$ 
(since $\pi^*V_3$ is a direct sum of stable bundles), 
hence there exists $d$ such that the map
\begin{equation}\label{mapE'}
\oplus_{x\in E'_d} H^0(E', L\otimes\PP_x)\otimes
\Hom(L\otimes\PP_x,\pi^*V_3)\ra H^0(E',\pi^*V_3)
\end{equation}
is surjective. Now we notice that for every $x\in E'$
there is a natural commutative diagram
\begin{equation}\label{comm}
\begin{array}{ccc}
H^0(E',L\otimes\PP_x)\otimes \Hom(L\otimes\PP_x,\pi^*V_3)&\lrar{} 
& H^0(E',\pi^*V_3)\\
\ldar{}&&\ldar{}\\
H^0(E,\pi_*(L\otimes\PP_x))\otimes \Hom(\pi_*(L\otimes\PP_x),V_3)&\lrar{}
& H^0(E,V_3)
\end{array}
\end{equation}
in which the right vertical arrow is the following composition of
natural morphisms:
$$H^0(E',\pi^*V_3)\wt{\ra} H^0(E,\pi_*\pi^*V_3)\ra H^0(E,V_3).$$
Since $V_3$ is a direct summand in $\pi_*\pi^*V_3$,
this map is surjective.
Also for every $y\in E$ we have an isomorphism
$$\pi_*(L\otimes\PP_{\hat{\pi}(x)})\simeq (\pi_* L)\otimes\PP_x$$
where $\hat{\pi}:E\ra E'$ is the isogeny dual to $\pi$.
Hence, the surjectivity of (\ref{mapE'}) implies the surjectivity
of the following map 
$$\oplus_{y\in \hat{\pi}^{-1}(E'_d)} H^0(E, (\pi_* L)\otimes\PP_y)\otimes
\Hom((\pi_*L)\otimes\PP_y,V_3)\ra H^0(E,V_3)$$
as required.

It remains to prove the statement for a
triple $(\O_E,L,V)$ where $L$ is a line bundle,
$V$ is a stable bundle, such that $0<\deg(L)<\mu(V)$. 
Note that $H^0(E,V)$ is an irreducible
representation of the Heisenberg group $H$ which is an extension of
$E_d$ by $\G_m$, where $d=\deg V$. More precisely, $H$ is the group of pairs
$(x,\phi)$ where $x\in E_d$, $\phi:V\ra t_x^*V$. It follows that  
the image of the natural map
$$\oplus_{x\in E_d} H^0(E,t_x^*L)\otimes\Hom(t_x^*L,V)\ra H^0(E,V)$$
is invariant under the $H$-action. Therefore, it suffices to prove that this
map is not zero. Let $f:L\ra V$ be a non-zero morphism.
Then it is an injection of sheaves, hence 
the induced morphism $H^0(E,L)\ra H^0(E,V)$ is injective which finishes
the proof.
\ed

Using the equivalence of the Fukaya category of $\R^2/\Z^2$ 
(without higher products)
with the category of vector bundles on $E$ we deduce
the following corollary.

\begin{cor}\label{surj}  
Let $(L,t)$ (resp. $(L',t')$) be an object in the
Fukaya category, where $L$ (resp. $L'$) is a line of slope $\la$
(resp. $\la'$). Assume that $\la<\mu<\la'$, where $\mu\in\Q$.
Then there exists a finite number of objects $(M_i,t_i)$ in
Fukaya category, where $M_i$ are lines of slope $\mu$, such that
the composition
$$m_2:\oplus_i\Hom((L,t),(M_i,t_i))\otimes
\Hom((M_i,t_i),(L',t'))\ra\Hom((L,t),(L',t'))$$ 
is surjective. Furthermore, one can choose $(M_i,t_i)$ in
a generic position (i.e.
in such a way that $M_i$ do not pass through a finite number of 
given points).
\end{cor}

\section{Expression of all triple products via trapezoid ones}

In this section we use the $A_{\infty}$-axiom (\ref{ainfty}) to get
a simple expression of an arbitrary triple product in $\FF_s$ (corresponding to a 
transversal configuration of lines)
in terms of trapezoid triple products and all double products. 

First, consider the triple product $m_3(r,s,t)$ where
$r\in \Hom^1(L_0,L_1)$, $s\in\Hom^0(L_1,L_2)$, $t\in\Hom^0(L_2,L_3)$,
$\la_0>\la_1<\la_2<\la_3$, $\la_0<\la_3$. If $\la_2=\la_0$ 
then this is a trapezoid product.
Otherwise, there are two possibilities:

\noindent a) $\la_2>\la_0$. In this case using Corollary \ref{surj} 
we can write $s$ as a linear combination of products $s's''$ where
$s'\in \Hom^0(L_1,L'_0)$, $s''\in \Hom^0(L'_0,L_2)$,
$L'_0$ is a line of slope $\la_0$, $L'_0\neq L_0\mod(\Z^2)$.
Then applying the 
$A_{\infty}$-constraint to the quadruple $r,s',s'',t$ we obtain 
$$m_3(r,s's'',t)= -m_3(r,s',s'')t + m_3(r,s',s''t)$$ 
(note that $rs'=0$ by assumption while
$m_3(s',s'',t)=0$ as an element of $\Hom^{-1}(L_1,L_3)$).

\noindent b) $\la_2<\la_0$. In this case we write $t$ as a linear
combination of products $t't''$ where
$t'\in \Hom^0(L_2,L'_0)$, $t''\in \Hom^0(L'_0,L_3)$,
$L'_0$ is a line of slope $\la_0$, $L'_0\neq L_0\mod(\Z^2)$.
Then we have $m_3(r,s,t')\in\Hom^0(L_0,L'_0)=0$. 
Hence, applying the $A_{\infty}$-constraint
to the quadruple $r,s,t',t''$ we get
$$m_3(r,s,t't'')=- m_3(rs,t',t'') + m_3(r,st',t'')$$ 
(note that $m_3(s,t',t'')=0$).

One deals similarly with products of the type
$$\Hom^0(L_0,L_1)\otimes\Hom^0(L_1,L_2)\otimes\Hom^1(L_2,L_3)\ra
\Hom^0(L_0,L_3)$$
where $\la_0<\la_1<\la_2>\la_3$.

Now let us consider $m_3(r,s,t)$ where
$r\in\Hom^0(L_0,L_1)$, $s\in\Hom^1(L_1,L_2)$,
$t\in\Hom^0(L_2,L_3)$, $\la_0<\la_1>\la_2<\la_3$, $\la_0<\la_3$.
If $\la_1=\la_3$ then this is a trapezoid product.
Otherwise, there are two possibilities:

\noindent 
a) $\la_1>\la_3$. Then we can write $r$ as a linear
combination of products $r'r''$, where $r'\in\Hom^0(L_0,L'_3)$,
$r''\in\Hom^0(L'_3,L_1)$, $L'_3$ is a line of slope $\la_3$,
$L'_3\neq L_3\mod(\Z^2)$. Applying $A_{\infty}$-constraint to
$(r',r'',s,t)$ we get
$$m_3(r'r'',s,t)= m_3(r',r''s,t) - m_3(r',r'',st) +
m_3(r',r'',s)t$$
(since $m_3(r'',s,t)\in\Hom^0(L'_3,L_3)=0$). The first two terms
in the RHS are trapezoid, while the product $m_3(r',r'',s)$ is of
the form considered before.

\noindent 
a) $\la_1<\la_3$. In this case we can write $t$ as a linear
combination of products $t't''$, where $t'\in\Hom^0(L_2,L'_1)$,
$t''\in\Hom^0(L'_1,L_3)$, $L'_1$ is a line of slope $\la_1$,
$L'_1\neq L_1\mod(\Z^2)$. Applying $A_{\infty}$-constraint to
$(r,s,t',t'')$ we get
$$m_3(r,s,t't'')= - m_3(rs,t',t'') + r m_3(s,t',t'') +
m_3(r,s,t')t''$$
(since $st'\in\Hom^1(L_1,L'_1)=0$). The products 
$m_3(s,t',t'')$ and $m_3(r,s,t')$
are trapezoid, while the product $m_3(rs,t',t)$ is of
the form considered before.

Finally, using 
the cyclic symmetry of $m_3$ we can reduce all non-zero
transversal higher products $m_3$ to the ones considered above.
For example, the product
$$\Hom^1(L_0,L_1)\otimes \Hom^1(L_1,L_2)\otimes \Hom^0(L_2,L_3)\ra 
\Hom^1(L_0,L_3)
$$
is equivalent to the product
$$\Hom^1(L_1,L_2)\otimes \Hom^0(L_2,L_3)\otimes \Hom^0(L_3,L_0)\ra
\Hom^0(L_1,L_0).$$

\bigskip

\noindent
{\it Proof of Theorem \ref{main}.}

\noindent
First, let us introduce some notation. We assume that for all objects
$(L,t)$ of the Fukaya category that appear below a representative
of the line $L$ modulo $\Z^2$-translations is fixed. Then for
every pair of objects $(L_1,t_1)$ and $(L_2,t_2)$ 
and every intersection point
$P\in (L_1+\Z^2)\cap (L_2+\Z^2)/\Z^2$ 
we represent $P$ in the
form $P_{12}(y_1,y_2+m\la_2+n)$, where $L_i=L(\la_i,y_i)$, and
set
$$e(P)=e_{12}(m,n)$$
where $e_{12}(m,n)$ is defined by formula (\ref{eij}).

As was explained in section \ref{triple} (after passing to a commensurable
lattice $\La$) we can assume that the lattice $\La$, the quadratic
form $Q$, and the cone $C$ come from a quadruple of rational numbers
$\la_0,\ldots,\la_3$ such that $\la_1<\la_0<\la_2<\la_3$.
Let us represent the variable $\zz\in\La_{\C}$ in the
form $\zz=\tau \vv-\ww$ with $\vv,\ww\in\La_{\R}$. Then
the value of $\Theta_{\La,Q,C}(\zz,\tau)$ appears
as the coefficient with $e(P_{03})$ of 
the triple product $m_3(e(P_{01}),e(P_{12}),e(P_{23}))$
for a quadruple of objects $(L_i,t_i)$, $i=0,\ldots,3$, 
and intersection points $P_{i,i+1}\in L_i\cap L_{i+1}$,
where $L_0=L(\la_0,0)$, $L_1=L(\la_1,0)$,
$L_2=(\la_2,(\la_2-\la_1)v_1)$, $L_3=(\la_3,(\la_0-\la_3)v_0)$,
$t_0=t_1=0$, $t_2=(\la_2-\la_1)w_1$, $t_3=(\la_0-\la_3)w_0$.
Now we can apply the above procedure of expressing this triple product 
in terms of the trapezoid ones. More precisely, due to the inequalities
$\la_1<\la_0<\la_2<\la_3$ we apply the very first case of the above
argument. This means that we choose a finite number of (not necessarily
distinct) objects
$(M_j,t_j)$, where $M_j$ are lines of slope $\la_0$
different from $L_0$, and intersection
points $Q_j\in L_1\cap (M_j+\Z^2)$, $R_j\in M_j\cap (L_2+\Z^2)$
such that $m_2(e(Q_j),e(R_j))$ form a basis in 
$\Hom((L_1,t_1),(L_2,t_2))$. 
The transition matrix from this basis to the standard basis
of intersection points of $L_1$ and $L_2$ modulo $\Z^2$
is given by elliptic functions of $(z_1,\tau)$. Thus, we can write
$$e(P_{12})=\sum_j \phi_j m_2(e(Q_j),e(R_j))$$
where $\phi_j$ are meromorphic elliptic
functions of $(z_1,\tau)$. Applying the 
$A_{\infty}$-identity (\ref{ainfty}) we get
\begin{align*}
&m_3(e(P_{01}),e(P_{12}),e(P_{23}))=\\
&\sum_j \rho_j(- m_3(e(P_{01}),e(Q_j),e(R_j))e(P_{23})+
m_3(e(P_{0,1}),e(Q_j),e(R_j)e(P_{23}))).
\end{align*}
Now by the results of section \ref{triple} and by formula
(\ref{appell})
the coefficients of the trapezoid product $m_3(e(P_{0,1}),e(Q_j),e(R_j))$  
are given (up to factors of the form $\exp(\pi i a\tau)$ with $a\in\Q_{>0}$)
by the functions of the form 
$\kappa_{e}(a\tau+b,gz_1+c\tau+d;h\tau)$, 
where $a,b,c,d,e,g\in\Q$, $h\in\Q_{\>0}$, 
Multiplying the result with $e(P_{2,3})$ means that we get
some linear combination of the above functions with coefficients
which are elliptic functions of $(s(\zz),\tau)$ for some linear functional $s$.
Finally the products
$m_3(e(P_{0,1}),e(Q_j),e(R_j)e(P_{2,3}))$
are expressed via elliptic functions of $(s(\zz),\tau)$ and
the functions of the form
$\kappa_e(a\tau+b,gz_2+c\tau+d;h\tau)$.
\ed

\begin{rems} 1. Since after rescaling $Q$ the slopes $\la_i$ can
always be chosen to be integers one can replace the reference
to Proposition \ref{surj-pr} in the above proof by the well-known
surjectivity statement for morphisms between line bundles.

\noindent
2. It may seem strange that in Theorem \ref{main} we
substitute only
constants in the second argument of $\kappa$. However, the function
$\kappa$ satisfies some identities (see \cite{H}, p. 481, formula (45)
and the next one, or \cite{P}, formula (3.4.3)) which imply that one
can express $\kappa(y,x;\tau)$ in terms of $\tau$-elliptic functions
and the function $\kappa(c,x+y-c;\tau)$ for any $c\in\Q+\Q\tau$.
\end{rems}

\section{Example}

In this section we will give an example of identity produced by
Theorem \ref{main}.
Let us fix $a\in\Z$ such that $a\ge 2$ and consider the
quadratic form $Q$ on $\Z^2$ given by
$$Q(n_0,n_1)=an_0^2+4an_0n_1+(4a-2)n_1^2.$$
Let us also consider the following split quadratic forms:
$$Q^1(n_0,n_1)=2(n_0+n_1)n_1$$
on $\Z^2$ and
$$Q^2(n_0,n_1)=(2n_0+\frac{2a-1}{a}n_1)n_1$$
on the lattice $\La^2=\{(n_0,n_1):\ n_0\in\Z, n_1\in a\Z\}$.
As a cone $C$ in all three cases we choose $x_0x_1>0$ 
(with $x_0>0$ in $C^+$) and set
$$\Theta(z_0,z_1)=\Theta_{\Z^2,Q,C}(z_0,z_1;\tau),$$
$$\Theta^1(z_0,z_1)=\Theta_{\Z^2,Q^1,C}(z_0,z_1;\tau),$$
$$\Theta^2_{c_0,c_1}(z_0,z_1)=\Theta_{\La^2,Q^2,C;(c_0,c_1)}
(z_0,z_1;\tau)$$
(we omit the variable $\tau$ in notation for brevity).
Let us also denote by $\De(z)$
the determinant of the $2\times 2$ matrix
$$(\theta_{2\Z,i}(z+\frac{j}{2},\frac{\tau}{2}))_{i\in\Z/2\Z,j\in\Z/2\Z}.$$
Then we have the following identity:
\begin{align*}
&\De(z_1-\frac{\tau}{2})\Theta(z_0,z_1)=
\theta_{2\Z,1}(z_1-\frac{\tau-1}{2},\frac{\tau}{2})\times\\
&\{-\theta_{a\Z}(z_0+2z_1,\frac{\tau}{a})\Theta^1(-2z_1,\frac{\tau}{2})+
\sum_{l\in\Z/a\Z}\theta_{a\Z,l}(-z_0-2z_1+\frac{\tau}{2a},\frac{\tau}{a})
\Theta^2_{\frac{a-1}{a}l,-l}(z_0,\frac{\tau}{2})\}\\
&-\theta_{2\Z,1}(z_1-\frac{\tau}{2},\frac{\tau}{2})\times\\
&\{-\theta_{a\Z}(z_0+2z_1,\frac{\tau}{a})\Theta^1(-2z_1,\frac{\tau-1}{2})+
\sum_{l\in\Z/a\Z}\theta_{a\Z,l}(-z_0-2z_1+\frac{\tau-1}{2a},\frac{\tau}{a})
\Theta^2_{\frac{a-1}{a}l,-l}(z_0,\frac{\tau-1}{2})\}.
\end{align*}

Let us write the variables in
the form $z_i=v_i\tau-w_i$, where $i=0,1$, and
consider the following objects in the Fukaya category:
$(L_0=L(0,0),t_0=0)$, $(L_1=L(-1,0),t_1=0)$,
$(L_2=L(1,2v_1),t_2=2w_1)$, $(L_3=L(\frac{a}{a-1},-\frac{a}{a-1}v_0),
t_3=-\frac{a}{a-1}w_0)$. Then
the series $\Theta(z_0,z_1)$ is equal to the coefficient with
$e_{03}$ in the triple product
$m_3(e_{01},e_{12},e_{23})$ where we set $e_{ij}=e_{ij}(0,0)$.
Now we consider two auxiliary objects in the Fukaya category:
$(M_0=L(0,\frac{1}{2}), 0)$ and $(M_1=L(0,\frac{1}{2}), \frac{1}{2}dx)$.
There are unique points of intersection
$Q_i\in L_1\cap M_i$, $R_i\in M_i\cap L_2$, $i=0,1$.
Note that the points $Q_0$ and $Q_1$ (resp. $R_0$ and $R_1$) coincide
but we denote them differently since they belong to morphism spaces between
different objects in the Fukaya category. The first step is
to represent $e_{12}$ as a linear combination of
$e(Q_0)e(R_0)$ and $e(Q_1)e(R_1)$. We have
$$e(Q_i)e(R_i)= \theta_{2\Z}(z_1-\frac{\tau-i}{2},\frac{\tau}{2})e_{12}+
\theta_{2\Z,1}(z_1-\frac{\tau}{2},\frac{\tau-i}{2})e_{12}(0,1)$$
for $i=0,1$. Hence,
$$\Delta(z_1-\frac{\tau}{2})e_{12}=
\theta_{2\Z,1}(z_1-\frac{\tau-1}{2},\frac{\tau}{2})e(Q_0)e(R_0)-
\theta_{2\Z,1}(z_1-\frac{\tau}{2},\frac{\tau}{2})e(Q_1)e(R_1).$$
Next we use the formula
$$m_3(e_{01},e(Q_i)e(R_i),e_{23})=-m_3(e_{01},e(Q_i),e(R_i))e_{23}+
m_3(e_{01},e(Q_i),e(R_i)e_{23}).$$
We have
$$m_3(e_{01},e(Q_i),e(R_i))=\Theta^1(-2z_1,\frac{\tau-i}{2})e_{02},$$
$$e(R_i)e_{23}=\sum_{l\in\Z/a\Z}
\theta_{a\Z,l}(-z_0-2z_1+\frac{\tau-i}{2a},\frac{\tau}{a})e^i_{03}(0,-l)$$
where $i=1,2$,
$e^i_{03}(0,l)$ are the basis elements in $\Hom(M_i,L_3)$ defined by 
(\ref{eij}). The coefficient with $e_{03}$ in the product
$e_{02}e_{03}$ is equal to
$$\theta_{a\Z}(z_0+2z_1,\frac{\tau}{a}).$$
One more computation shows that the coefficient
with $e_{03}$ in the product $m_3(e_{01},e(Q_i),e^i_{03}(0,-l))$
(where $i=1,2$) is equal to
$$\Theta^2_{\frac{a-1}{a}l,-l}(z_0,\frac{\tau-i}{2}).$$
Combining all these calculations we get the identity above.

\section{Massey products} 

In this section we consider a family of well-defined univalued triple Massey
products which define global sections of certain line bundles on the
second cartesian power of the universal curve over the moduli stack of 
elliptic curves with some level structure. 

\subsection{Definition of triple Massey products}
\label{def-mas}

Let $V_i$, $0\le i\le 3$ be holomorphic vector bundles on a complex manifold.
Let $\a_1\in\Hom(V_0,V_1)$, $\a_2\in\Ext^1(V_1,V_2)$, $\a_3\in\Hom(V_2,V_3)$
be elements satisfying $\a_2\circ\a_1=0$, $\a_3\circ\a_2=0$.
Below we recall two equivalent constructions of the triple Massey product
$MP(\a_1,\a_2,\a_3)$ which belongs to the cokernel of the morphism
\begin{equation}\label{masmor}
\Hom(V_0,V_2)\oplus\Hom(V_1,V_3)\ra\Hom(V_0,V_3):(\b_1,\b_2)\mapsto 
\a_3\circ\b_1+\b_2\circ\a_1.
\end{equation}

Let us represent $\a_2$ by a $\dbar$-closed $(0,1)$-form $\wt{\a}_2$ with 
values 
in $V_1^*\otimes V_2$. Then by our assumption we have
$$\wt{\a}_2\circ\a_1=\dbar(\a_{12}),$$
$$\a_3\circ\wt{\a}_2=\dbar(\a_{23})$$
for some sections $\a_{12}\in C^{\infty}(V_1^*\otimes V_2)$,
$\a_{23}\in C^{\infty}(V_2^*\otimes V_3)$.
Now we set
$$MP(\a_1,\a_2,\a_3)=\a_3\circ\a_{12}-\a_{23}\circ\a_1.$$
The ambiguity in a choice of $\a_{12}$ and $\a_{23}$ precisely means
that $MP(\a_1,\a_2,\a_3)$ is correctly defined modulo the image of the map 
(\ref{masmor}).

In the second definition 
\footnote{This definition is a particular case of the general construction
of Massey products in triangulated categories, cf. \cite{GM},IV.2}
we consider an extension
$$0\ra V_2\stackrel{i}{\ra} V\stackrel{p}{\ra} V_1\ra 0$$ 
with the class $\a_2$. By our assumption there exist morphisms
$\a'_1:V_0\ra V$ and $\a'_3:V\ra V_3$ such that
$$\a_1=p\circ\a'_1,$$
$$\a_3=\a'_3\circ i.$$
Now the composition $\a'_3\circ\a'_1\in\Hom(V_0,V_3)$ is well-defined
modulo the image of (\ref{masmor}).

\begin{prop}\label{masprop} One has
$$MP(\a_1,\a_2,\a_3)=-\a'_3\circ\a'_1$$
in the cokernel of the map (\ref{masmor}).
\end{prop}

\Pf . A choice of a closed $(0,1)$-form $\wt{\a}_2$ representing $\a_2$
leads to the choice of $V$ as follows. We set $V=V_2\oplus V_1$ as a
$C^{\infty}$-bundle and define a holomorphic structure on it by
the following $\dbar$-operator:
$$\dbar_V=\left(\matrix \dbar_{V_2} & \wt{\a}_2 \\ 0 & \dbar_{V_1}
\endmatrix \right).$$
Let $\si:V_1\ra V$ be the natural $C^{\infty}$-splitting arising from this 
description (so $p\circ\si=\id_{V_1}$). Then $\dbar(\si)=i\circ\wt{\a}_2$.
Define $\rho:V\ra V_2$ by the condition 
$$i\circ\rho=\si\circ p-\id_V.$$
Then we have
$$\dbar(\rho)=\wt{\a}_2\circ p.$$
Thus, we can choose 
$$\a_{12}=\rho\circ\a'_1,$$
$$\a_{23}=\a'_3\circ\si.$$
Hence,
$$\a_3\circ\a_{12}-\a_{23}\circ\a_1=\a'_3\circ i\circ\rho\circ\a'_1-
\a'_3\circ\si\circ p\circ\a'_1=-\a'_3\circ\a'_1.$$
\ed

We will be interested in a particular case when 
$\Hom(V_0,V_2)=\Hom(V_1,V_3)=0$.
In this case the Massey product $MP(\a_1,\a_2,\a_3)$ is an element of
$\Hom(V_0,V_3)$.
Homological mirror conjecture for elliptic curve $E=E_{\tau}$
(proven in \cite{P-hmc} for transversal products) implies that
$MP(\a_1,\a_2,\a_3)$ is equal to the corresponding
triple Fukaya product $m_3(\a_1,\a_2,\a_3)$. Indeed, the corresponding products 
are homotopic
but in our case the homotopy takes values in zero spaces (see \cite{P-YB}, 
sec.1.1, for a more detailed explanation).

\begin{rem} The notion of transversality considered in \cite{P-hmc} has
to be strengthened, since the definition of Fukaya products $m_k$ given in
section \ref{Fukdef-sec}
requires that no three of the corresponding circles intersect in one point (this 
was overlooked in \cite{P-hmc}). However, all the proofs of \cite{P-hmc}
can be easily modified accordingly.
\end{rem}

\subsection{Massey products for line bundles}
\label{masline-sec}

Let us fix a quadruple of integers $(d_0,d_1,d_2,d)$ such that
$1\le d_0<\min(d_1,d_2)$, $d_0+d=d_1+d_2$. 
Let $\LL_0$, $\LL_1$, $\LL_2$ and $\LL$ 
be line bundles on elliptic curve $E$ of degrees
$d_0$, $d_1$, $d_2$ and $d$ respectively, such that $\LL_0\LL\simeq
\LL_1\LL_2$ (here and below we skip the sign of tensor product
between line bundles for brevity).
Then for any pair of sections $s_1\in H^0(E,\LL_1)$, $s_2\in H^0(E,\LL_2)$
and an element $e\in H^1(\LL^{-1})=H^1(E,\LL_0\LL_1^{-1}\LL_2^{-1})$
such that the compositions
$s_1e\in H^1(E,\LL_0\LL_2^{-1})$ and $es_2\in H^1(E,\LL_0\LL_1^{-1})$
are zero, the triple Massey product $MP(s_1,e,s_2)$  
defined in \ref{def-mas}
is an element of $H^0(E,\LL_0)$ 
(since $H^0(E,\LL_0\LL_1^{-1})=H^0(E,\LL_0\LL_2^{-1})=0$).
Thus, if we denote by $K_{s_1,s_2}\subset H^1(E,\LL^{-1})$ 
the kernel of the natural map
$$H^1(E,\LL^{-1})\ra
H^1(E,\LL_0\LL_1^{-1})\oplus H^1(E,\LL_0\LL_2^{-1}):e\mapsto(s_1e,es_2)$$
then the Massey product defines a linear map
\begin{equation}\label{MPK}
MP:K_{s_1,s_2}\ra H^0(E,\LL_0).
\end{equation}
Assume that $s_1\neq 0$, $s_2\neq 0$.
Let $D$ be the divisor of common zeroes of $s_1$ and $s_2$.
Let us denote $\LL'_i=\LL_i(-D)$ for $i=0,1,2$, $\LL'=\LL(-D)$,
so that we still have $\LL'_0\LL'\simeq\LL'_1\LL'_2$.
Then we have the induced sections $s'_1\in H^0(E,\LL'_1)$ and
$s'_2\in H^0(E,\LL_2')$ and a canonical map
$$\varphi:K_{s_1,s_2}\ra K_{s'_1,s'_2}$$
induced by the morphism $\LL^{-1}\ra\LL^{-1}(D)=(\LL')^{-1}$.
On the other hand, we can consider $H^0(E,\LL'_0)=H^0(E,\LL_0(-D))$
as a subspace in $H^0(E,\LL_0)$.

\begin{lem}\label{divisor}
For any $e\in K_{s_1,s_2}$ one has
$$MP(s_1,e,s_2)=MP(s'_1,\varphi(e),s'_2).$$
\end{lem}

\Pf . First let us notice that
$$H^0(E,\LL'_0(\LL'_i)^{-1})=H^0(E,\LL_0\LL_i^{-1})=0$$
for $i=1,2$, so $MP(s'_1,\varphi(e),s'_2)$ is defined.
Let
$$0\ra \LL^{-1}\ra E\ra \O_E\ra 0$$
be an extension representing $e$. Let $f:\LL_1^{-1}\ra E$
be the lifting of $s_1:\LL_1^{-1}\ra\O_E$. Then the image of
$f$ belongs to the following subextension $E'$:
$$0\ra \LL^{-1}\ra E'\ra \O_E(-D)\ra 0.$$
Note that this extension represents the class $\varphi(e)$.
Let $g':E'\ra\LL_0\LL_1^{-1}(-D)$ be the lifting of the
map $s'_2:\LL^{-1}\ra\LL_0\LL_1^{-1}(-D)$.
Then according to Proposition \ref{masprop} we have
$$MP(s'_1,\varphi(e),s'_2)=-g'\circ f.$$
On the other hand, the push-out of the extension $E$ by
the morphism $\LL^{-1}\ra\LL^{-1}(D)$
coincides with extension
$$0\ra \LL^{-1}(D)\ra E'(D)\ra \O_E\ra 0.$$
Thus, we have an embedding $i:E\ra E'(D)$, such that
$i|_E':E'\ra E'(D)$ is the natural map. In particular,
we can take $g=g'\circ i:E\ra\LL_0\LL_1^{-1}$ as the lifting
of the map $s_2:\LL^{-1}\ra\LL_0\LL_1^{-1}$. 
Applying Proposition \ref{masprop} again we obtain
$$MP(s_1,e,s_2)=-g\circ f$$
which finishes the proof.
\ed

Thus, it suffices to study the case when
the sections $s_1$ and $s_2$ have no common zeroes.
In this case one has an exact sequence
\begin{equation}\label{extension}
0\ra \LL_1^{-1}\LL_2^{-1}\stackrel{\a}\ra \LL_1^{-1}\oplus\LL_2^{-1}
\stackrel{\b}{\ra} \O_E\ra 0
\end{equation}
where $\a(s)=(-ss_2,ss_1)$, $\b(t_1,t_2)=t_1s_1+t_2s_2$.
Tensoring by $\LL_0$
and considering the corresponding sequence of cohomologies we get
the following exact complex $C$:
$$0\ra H^0(E,\LL_0)\stackrel{\delta_{-1}}{\ra} H^1(E,\LL^{-1})
\stackrel{\delta_0}{\ra}
H^1(E,\LL_0\LL_1^{-1})\oplus H^1(E,\LL_0\LL_2^{-1})\ra 0$$

\begin{prop}\label{mas}
For any $t\in H^0(E,\LL_0)$ one has
$$MP(s_1,\delta_{-1}(t),s_2)=t.$$
\end{prop}

\Pf . The class $\delta_{-1}(t)\in\Ext^1(\LL_0^{-1},\LL_1^{-1}\LL_2^{-1})$
is represented by the extension
$$0\ra\LL_1^{-1}\LL_2^{-1}\ra E\ra \LL_0^{-1}\ra 0$$
where $E=\b^{-1}(t(\LL_0^{-1})).$
The morphism $\LL_0^{-1}\LL_1^{-1}\stackrel{s_1}{\ra}\LL_0^{-1}$
lifts to the morphism 
$$f:\LL_0^{-1}\LL_1^{-1}\ra E:s\mapsto (st,0),$$
while the morphism
$\LL_1^{-1}\LL_2^{-1}\stackrel{s_2}{\ra}\LL_1^{-1}$
lifts to the morphism
$$g:E\ra\LL_1^{-1}:(t_1,t_2)\mapsto -t_1.$$
According to Proposition \ref{masprop} we have
$$MP(s_1,\delta_{-1}(t),s_2)=-g\circ f,$$
so the result follows from the above formulas for $f$ and $g$.
\ed

The spaces $K_{s_1,s_2}$ can be considered as stalks of the sheaf $\KK$ 
over the appropriate moduli stack. The Massey product gives a morphism
from $\KK$ to the bundle with the fibre $H^0(E,\LL_0)$. The previous
proposition shows that over an open part this is an isomorphism inverse
to $\delta_{-1}$. 
In the section \ref{masgen} we will compute the map $MP$ in terms of indefinite
theta series using (\ref{m3}) and applying the above proposition we will get
a proof of Theorem \ref{linsys}.
In order to get modular (or Jacobi) forms from the coefficients of the map $MP$
we can try to map various standard line bundles (trivialized by theta functions)
to $\KK$.  
Below we will present two ways to do it. The first uses the determinants
and gives Jacobi forms.
The second approach (see section \ref{ex-mod}) is more direct and produces
modular forms but it requires additional assumptions on the integers 
$(d_0,d_1,d_2)$.

\subsection{Determinantal approach}\label{det-Mas-sec}
  
Let us consider the $1$-dimensional vector space
$$M=\det H^1(E,\LL^{-1})\otimes
\det H^1(E,\LL_0\LL_1^{-1})^*\otimes\det H^1(E,\LL_0\LL_2^{-1})^*$$
where for a vector space $V$ we denote by $\det V$ its top-degree wedge power.
We claim that for fixed $s_1$ and $s_2$ there is a canonical map
$$e_{s_1,s_2}:M\otimes\sideset{}{^{d_0-1}}\We
H^1(E,\LL^{-1})^*\ra K_{s_1,s_2}$$
Indeed, in general for a linear map $f:V\ra W$
between vector spaces of dimensions
$n$ and $n-k$ one can construct a linear map
$$k_f:\det V\otimes\det W^*\otimes\sideset{}{^{k-1}}\We V^*\ra V$$
such that its image belongs to $\ker(\phi)$ as follows.
Start with the morphism
$$\sideset{}{^{n-k}}\We f^*:\det W^*=\sideset{}{^{n-k}}\We W^*\ra
\sideset{}{^{n-k}}\We V^*$$
and then consider the following composition
$$\phi:\det W^*\otimes\sideset{}{^{k-1}}\We V^*\ra
\sideset{}{^{n-k}}\We V^*\otimes\sideset{}{^{k-1}}\We V^*\ra
\sideset{}{^{n-1}}\We V^*$$
where the second arrow is given by the wedge product. It remains to
use the isomorphism $\We^{n-1} V^*\simeq V\otimes \det V^*$.
It is easy to see that the image of $k_f$ belongs to $\ker(f)$.
If $f$ is not surjective then $k_f=0$, otherwise, $k_f$ surjects
onto $\ker(f)$.
The map $e_{s_1,s_2}$ is defined by applying this construction to the map
$$
\delta_0:H^1(E,\LL^{-1})\ra H^1(E,\LL_0\LL_2^{-1})\oplus
H^1(E,\LL_0\LL_1^{-1})
$$
induced by compositions with $s_1$ and $s_2$
(note that dimensions of the spaces are $d_1+d+2-d_0$ and $d_1+d_2-2d_0$).
Composing $e_{s_1,s_2}$ with the map (\ref{MPK}) above we obtain a
linear map 
$$MP_{s_1,s_2}^{\det}:M\otimes
\sideset{}{^{d_0-1}}\We H^1(E,\LL^{-1})^*
\ra H^0(E,\LL_0):\xi\mapsto MP(s_1,e_{s_1,s_2}(\xi),s_2).$$
It is easy to see that if the sections $s_1$ and $s_2$ have a common zero
then $\delta_0$ is not surjective, hence, $MP_{s_1,s_2}^{\det}=0$. 

We refer to \cite{GKZ}, Appendix A, for the definition of
determinants of complexes.

\begin{prop}\label{mas2}
Assume that the divisors of $s_1$ and $s_2$ do not intersect. 
Then the map $\pm MP_{s_1,s_2}^{\det}$ is equal to the composition
$$M\otimes
\sideset{}{^{d_0-1}}\We H^1(E,\LL^{-1})^*
\ra M\otimes \sideset{}{^{d_0-1}}\We H^0(E,\LL_0)^*\ra H^0(E,\LL_0)$$
where the first arrow is induced by $\delta_{-1}^*$ and the second
arrow is induced by the isomorphism
$M\simeq \det H^0(E,\LL_0)$ given by the canonical trivialization
of the determinant of the exact complex $C$.
\end{prop}

\Pf . This follows from Proposition \ref{mas} and from the following
observation. For an exact complex of the form
$$0\ra K\stackrel{\delta_{-1}}{\ra} V\stackrel{\delta_0}{\ra} W\ra 0$$
the map
$$k_{\delta_0}:\det V\otimes\det W^*\otimes
\sideset{}{^{\dim K-1}}\We V^*\ra K$$
constructed above coincides
(up to a sign) with the composition
$$\det V\otimes\det W^*\otimes
\sideset{}{^{\dim K-1}}\We V^*\stackrel{\delta_0^*}{\ra}
\det V\otimes\det W^*\otimes
\sideset{}{^{\dim K-1}}\We K^*\ra K$$
where $\det V\otimes\det W^*\simeq\det K$ by the
canonical trivialization of the determinant of this
complex.
\ed

\begin{cor} Assume that $d_0=1$.
Fix some bases in the spaces $H^1(E,\LL^{-1})$,
$H^1(E,\LL_0\LL_i^{-1})$ and $H^0(E,\LL_0)$. Then we have
$$MP_{s_1,s_2}^{\det}=\pm\det(C,B)^{-1}$$
where $B$ is the corresponding basis of the complex $C$.
\end{cor} 

\subsection{Serre duality in homological mirror symmetry}

The identification of $\Ext^1$-spaces in the Fukaya category and the
category of bundles on elliptic curve $E=E_{\tau}$ uses the Serre duality
in the following form (see \cite{P-hmc}):
$$\Hom(V_1,V_2)\otimes\Ext^1(V_2,V_1)\ra\C:
f\otimes (gd\ov{z})\mapsto \int_E dz\wedge \Tr(f\circ gd\ov{z}).$$ 
In other words, this is a composition of the natural
map
$$\Hom(V_1,V_2)\otimes\Ext^1(V_2,V_1)\ra H^1(E,\O_E)$$
and the map
$$\phi:H^1(E,\O_E)\ra\C:\a\mapsto \int_E dz\wedge\a$$
where $\a$ is a $(0,1)$-form.
The map $\phi$ is in turn equal to the composition of the isomorphism
$$H^1(E,\O_E)\ra H^1(E,\om_E)$$
induced by the holomorphic $1$-form $\om_0=2\pi i dz$ and the functional
$$I:H^1(E,\om_E)\ra\C:\eta\mapsto-\frac{1}{2\pi i}\int_E\eta$$
where $\eta$ is a $(1,1)$-form.
The factor of $\frac{1}{2\pi i}$ in the definition of $I$ is important because 
then 
$I$
admits an algebraic definition (in particular, it is defined over the 
field of definition of $E$). In fact, we can define the functional 
$I:H^1(C,\om_C)\ra\C$ for any Riemann surface $C$ by the same formula. Then the 
following property of $I$ shows that it is algebraically defined.

\begin{lem} For every point $p\in C$ let $e_p\in H^1(C,\om_C)$ be
the class defined by the boundary homomorphism
$\C\simeq H^0(C,\O_p)\ra H^1(C,\om_C)$ coming from the
exact sequence
$$0\ra\om_C\ra\om_C(p)\stackrel{\Res_p}{\ra}\O_p\ra 0.$$
Then $I(e_p)=1$.
\end{lem}

\Pf . Let $U\subset C$ be an open disk containing $p$. 
Consider the Cech complex of $\om_C$ associated with
the covering $(U, C-p)$:
$$\CC^{\cdot}: \om_C(U)\oplus \om_C(C-p)\ra \om_C(U-p)$$
where the differential sends $(\a_U,\a_p)$ to $\a_U|_{U-p}-\a_p|_{U-p}$.
Since $H^1(C-p,\om_C)=H^1(U,\om_C)=0$ the complex $\CC^{\cdot}$ computes
the cohomology of $\om_C$. The residue map
$$\Res_p:\CC^1=\om_C(U-p)\ra\C$$
descends to the functional $I':H^1(C,\om_C)\ra\C$.
It suffices to prove that $I=I'$. To this end let us consider the Cech complex
associated with the same covering and with the complex of sheaves
$\Om^{1,0}\stackrel{\dbar}{\ra}\Om^{1,1}$ on $C$:
$$\CC\DD^{\cdot}:\Om^{1,0}(U)\oplus\Om^{1,0}(C-p)\stackrel{d_1}{\ra}
\Om^{1,1}(U)\oplus\Om^{1,1}(C-p)\oplus\Om^{1,0}(U-p)\stackrel{d_2}{\ra}
\Om^{1,1}(U-p)$$
where 
$$d_1(\a_U,\a_p)=(\dbar\a_U,\dbar\a_p,\a_U|_{U-p}-\a_p|_{U-p}),$$
$$d_2(\eta_U,\eta_p,\b)=\eta_U|_{U-p}-\eta_p|_{U-p}-\dbar\b.$$
The complex $\CC\DD^{\cdot}$ is concentrated in degrees $[0,2]$. We have
natural morphisms of complexes
$\CC^{\cdot}\ra\CC\DD^{\cdot}$ and $\Om^{1,\cdot}(C)\ra\DD^{\cdot}$ 
inducing isomorphisms on cohomologies.
Now we define the functional
$\wt{I}:\CC\DD^1\ra\C$ by the formula
$$\wt{I}(\eta_U,\eta_p,\b)=-\frac{1}{2\pi i}\int_{D}\eta_U-\frac{1}{2\pi i}
\int_{C-D}\eta_p+\frac{1}{2\pi i}\int_{\partial D}\b$$
where $D\subset U$ is a smaller disk containing $p$. 
It is easy to check that $\wt{I}\circ d_1=0$, hence,
$\wt{I}$ descends to a functional 
$$\wt{I}:H^1(C,\om_C)\simeq H^1(\CC\DD^{\cdot})\ra\C.$$
If $\eta$ is a global $(1,1)$-form then
$$\wt{I}(\eta|_U,\eta|_{C-p},0)=-\frac{1}{2\pi i}\int_C\eta,$$
hence $\wt{I}=I$. On the other hand, if $\a$ is a holomorphic $1$-form
on $U-p$ then
$$\wt{I}(0,0,\a)=\Res_p(\a),$$
hence $\wt{I}=I'$.
\ed

\subsection{Boundary homomorphism}\label{boundhom-sec}

Below we will need to calculate 
the matrices of $\delta_{-1}$ and $\delta_0$ with respect to
the standard bases of the terms of the complex $C$
(coming from Serre duality and the bases of theta functions in the spaces
of global sections) for any pair of sections $s_1,s_2$ with no common zeroes.
The components of $\delta_0$ are just the compositions 
with $s_1$ and $s_2$ so they are given by theta functions.
On the other hand, the map 
$$\delta_{-1}:H^0(E,\LL_0)\ra H^1(E,\LL^{-1})$$ 
is the composition
with the class $\delta(s_1,s_2)\in H^1(E,\LL_1^{-1}\LL_2^{-1})$ 
corresponding to the
extension (\ref{extension}). Via Serre duality $\delta_{-1}$ corresponds
to a bilinear form
$$B_{s_1,s_2}:H^0(E,\LL_0)\otimes H^0(E,\LL)\ra
H^0(E,\LL_1\LL_2)\ra\C$$
induced by the functional on $H^0(E,\LL_1\LL_2)$ dual to the class
$\delta(s_1,s_2)$
(recall that we always use the trivialization of $\om_E$ given by the
form $2\pi i dz$).

\begin{lem}\label{res} Assume that $s_1$ and $s_2$ have no common zeroes.
Then one has
$$B_{s_1,s_2}(s,t)=
\sum_{x\in Z(s_1)}\Res_{x}(\frac{2\pi i s(z)t(z) dz}{s_1(z)s_2(z)})$$
where $Z(s_1)$ is the divisor of zeroes of $s_1$.
\end{lem}

\Pf .
Applying the octahedron axiom to the composition
of arrows $\LL_1^{-1}\ra \LL_1^{-1}\oplus \LL_2^{-1}\stackrel{\b}{\ra}\O_E$ 
one can easily derive that the class $\delta(s_1,s_2)$ can
be represented as the following composition:
$$\O_E\ra \O_{Z(s_1)}\ra\LL_1^{-1}\LL_2^{-1}[1]$$
where the first arrow is the canonical one, the second arrow comes from
the exact triangle
$$\LL_1^{-1}\LL_2^{-1}\stackrel{s_1}{\ra}\LL_2^{-1}\ra
\O_{Z(s_1)}\ra\LL_1^{-1}\LL_2^{-1}[1]$$
(here we use the natural trivialization of $\LL_2|_{Z(s_1)}$ induced
by $s_2$).
In other words, this element in 
$\Ext^1(\O_{Z(s_1)},\LL_1^{-1}\LL_2^{-1})$
corresponds to the functional on
$H^0(E,\LL_1\LL_2|_{Z(s_1)})$
which at the point $x\in Z(s_1)$ is equal to 
$$\Res_{x}(\frac{2\pi i dz}{s_1(z)s_2(z)})$$
\ed

\section{Calculations}

In this section we will use homological mirror symmetry to relate the Massey
products considered in section \ref{masline-sec} to indefinite theta series.
We keep the notations of the previous section. 

\subsection{Proof of Theorem \ref{linsys}}
\label{masgen}

Let us fix a pair of complex numbers
$v_1\tau-w_1$ and $v_2\tau-w_2$
where $v_i,w_i\in\R$.
Now consider the following $4$ objects in the Fukaya category:
$(L_0=L(0,0),0)$, $(L_1=L(d_1,d_1v_1),d_1w_1)$, 
$(L_2=L(d_0-d_2,d_2v_2),d_2w_2)$, $(L_3=L(d_0,0),0)$.
Let $\LL_{\tau}$ be the basic line bundle on $E=E_{\tau}$
such that $\theta(z)=\theta(z,\tau)$ is a section of $\LL_{\tau}$. 
Under the equivalence with the category of bundles on $E$ our $4$
objects correspond to $\O_E$, $\LL_1$, $\LL_0\LL_2^{-1}$ and $\LL_0$
respectively, where
$\LL_1=t_{v_1\tau-w_1}^*\LL^{\otimes d_1}_{\tau}$, 
$\LL_2=t_{-v_2\tau+w_2}^*\LL^{\otimes d_2}_{\tau}$,
$\LL_0=\LL^{\otimes d_0}_{\tau}$.
Let $\La^+=\La^+(0,d_1,d_0-d_2,d_0)$ be the corresponding rank $2$ lattice.
Let us denote by $\xx=(x_0,x_1,x_2,x_3)$ the element of $\La^+\otimes\C$ such
that $x_0=v_1\tau-w_1$ and $x_3=v_2\tau-w_2$
(thus, $x_1=-\frac{(d-d_1)x_0+d_2x_3}{d}$, $x_2=-\frac{(d-d_2)x_3+d_1x_0}{d}$).

We start by fixing bases in all the relevant vector spaces.
According to section \ref{double} we have a natural isomorphism
$$\Hom(L_0,L_1)\simeq H^0(E,\LL_1)$$
which identifies the basis $(e_{01}(0,a), a\in\Z/d_1\Z)$
with $(\th_{d_1\Z,a}(z+x_0,\frac{\tau}{d_1}), a\in\Z/d_1\Z)$.
Similarly, 
$$\Hom(L_2,L_3)\simeq H^0(E,\LL_2)$$
such that $e_{23}(0,a)$, $a\in\Z/d_2\Z$, corresponds to
$\th_{d_2\Z,a}(z-x_3,\frac{\tau}{d_2})$.
On the other hand,
$$\Hom^1(L_1,L_2)\simeq H^1(E,\LL^{-1})\simeq
H^0(E,\LL)^*$$
in such a way that the basis
$(e_{12}(0,a),a\in\Z/d\Z)$ is dual to the basis
$(\th_{d\Z,-a}(z+x_0+x_1,\frac{\tau}{d}),a\in\Z/d\Z)$.
Finally, the space $\Hom^1(L_0,L_2)\simeq H^1(E,\LL_0\LL_2^{-1})$
(resp. $\Hom^1(L_1,L_3)\simeq H^1(E,\LL_0\LL_1^{-1})$) has the basis
$(e_{02}(0,a), a\in\Z/(d_2-d_0)\Z)$
(resp. $(e_{13}(0,a), a\in\Z/(d_1-d_0)\Z)$) and the space
$\Hom(L_0,L_3)\simeq H^0(E,\LL_0)$ has the basis
$(e_{03}(0,a), a\in\Z/d_0\Z)$ identified with 
$(\th_{d_0\Z,a}(z,\frac{\tau}{d_0}), a\in\Z/d_0\Z)$.

We are going to compute explicitly the map (\ref{MPK}) for
$s_1=\th_{d_1\Z}(z+x_0,\frac{\tau}{d_1})\in H^0(E,\LL_1)$
and $s_2=\th_{d_2\Z}(z-x_3,\frac{\tau}{d_2})\in H^0(E,\LL_2)$
(the corresponding morphisms in the Fukaya category are
$e_{01}=e_{01}(0,0)$ and $e_{23}=e_{23}(0,0)$).

First, using the formula (\ref{m2e}) one can easily compute the matrix of 
$\delta_0$.
Namely, we have
\begin{align*}
&m_2(e_{01},e_{12}(0,k))=
\sum_{n\in\Z/I_1}
\theta_{I_1,-\frac{k}{d}-n}(p_1x_1,p_1\tau)
e_{02}(0,k+dn),\\
&m_2(e_{12}(0,k),e_{23})=
\sum_{n\in\Z/I_2}\theta_{I_2,\frac{k}{d}+n}(p_2x_2,p_2\tau)
e_{13}(0,k+dn),
\end{align*}
where 
$I_i=\Z\cap\frac{d-d_i}{d_i}\Z$, $p_i=\frac{d_id}{d-d_i}$, $i=1,2$.
In these formulas we fix a representative of $k$ in $\Z$. 
The change of a representative
corresponds to a shift of the summation variable $n$. Here is a better way to 
write
these formulas which doesn't require a choice of a representative of $k$:
\begin{eqnarray}\label{comp}
m_2(e_{01},e_{12}(0,k))=
\sum_{n\in\Z/dI_1,\ n\equiv k(d)}
\theta_{I_1,-\frac{n}{d}}(p_1x_1,p_1\tau)
e_{02}(0,n),\nonumber\\
m_2(e_{12}(0,k),e_{23})=
\sum_{n\in\Z/dI_2,\ n\equiv k(d)}
\theta_{I_2,\frac{n}{d}}(p_2x_2,p_2\tau)
e_{13}(0,n),
\end{eqnarray}
Note that the coefficient of $m_2(e_{01},e_{12}(0,k))$ 
(resp. $m_2(e_{12}(0,k),e_{23})$) with
$e_{02}(0,i)$ (resp. $e_{13}(0,j)$) is non-zero only if $k\equiv i(d,d-d_1)$
 (resp. $k\equiv j(d,d-d_2)$).
Now it is easy to deduce that an element 
$e=\sum_{k\in\Z/d\Z}c_k e_{12}(0,k)$ (where $c_k\in\C$) 
belongs to the subspace $K_{s_1,s_2}=\ker(\delta_0)$ if and only
if
\begin{eqnarray}\label{conds}
\sum_{n\in\Z/dI_1,\ n\equiv i(d-d_1)}
c_n\theta_{I_1,-\frac{n}{d}}(p_1x_1,p_1\tau)=0,\nonumber\\
\sum_{n\in\Z/dI_2,\ n\equiv j(d-d_2)}
c_n\theta_{I_2,\frac{n}{d}}(p_2x_2,p_2\tau)=0  
\end{eqnarray}
for all $i\in\Z/(d-d_1)\Z$ and $j\in\Z/(d-d_2)\Z$.

Assume first that $x_1$ and $x_2$ are sufficiently generic, so that the 
corresponding circles in $\R^2/\Z^2$ form a transversal configuration.
Then we can use \eqref{m3} to compute the relevant triple Fukaya products.
Note that the projection $(n_0,n_1,n_2,n_3)\mapsto (n_0,n_3)$ defines
an isomorphism of the lattice $\La^+$ with the lattice
$$\La_{d_1,d_2,d}=\{ (m,n)\in\Z^2:\ d_1m\equiv d_2n (d)\}.$$
On the other hand, the projection $(n_0,n_1,n_2,n_3)\mapsto (n_1,n_2)$
maps $\La^+$ onto $\La_{d_1,d_2,d_0}$.
Thus, we have
\begin{align*}
&m_3(e_{01},e_{12}(0,k),e_{23})=\\
&\sum_{(n_1,n_2)\in \Z^2/\La_{d_1,d_2,d_0}}\Theta_{\La_{d_1,d_2,d},Q,mn>0;
(-\frac{k+d_2n_2-d_1n_1}{d_0}-n_1,
\frac{k+d_2n_2-d_1n_1}{d_0}-n_2)}(x_0,x_3;\tau)\\
&e_{03}(0,k+d_2n_2-d_1n_1)
\end{align*}
for $k\in\Z/d\Z$, where the form $Q$ is given by
$$Q(m,n)=\frac{1}{d}(d_1(d-d_1)m^2+2d_1d_2mn+(d-d_2)d_2n^2)$$

Let $(e_{03}^*(0,l), 0\le l<d_0)$ be the dual basis to $(e_{03}(0,l))$.
Then for $e=\sum_k c_k e_{12}(0,k)$ we have
$$\wt{F}_l:=\lan e_{03}^*(0,l),MP(s_1,e,s_2)\ran=\sum_{k} c_k
\lan e_{03}^*(0,l),m_3(e_{01},e_{12}(0,k),e_{23})\ran=
\sum_{k} C_{l,k} c_k$$
where $C_{l,k}$ is zero unless
$k\equiv l(d_0,d_1,d_2)$ in which case
$$C_{l,k}=\Theta_{\La_{d_1,d_2,d},Q,mn>0;(-\frac{l}{d_0}-i-n_1,
\frac{l}{d_0}+i-n_2)}(x_0,x_3;\tau)$$
where $n_1$, $n_2$, and $i$ are some integers satisfying
$$k+d_2n_2-d_1n_1=l+d_0 i.$$
Denoting $m_1=-i-n_1$, $m_2=i-n_2$ we can rewrite this formula as follows:
\begin{equation}\label{Clk}
C_{l,k}=\Theta_{\La_{d_1,d_2,d},Q,mn>0;(-
\frac{l}{d_0}+m_1,\frac{l}{d_0}+m_2)}(x_0,x_3;\tau)
\end{equation}
where $(m_1,m_2)\in\Z^2/\La_{d_1,d_2,d}$ is defined by the congruence
$$d_2m_2-d_1m_1\equiv k-l(d).$$
Thus, we have
\begin{align*}
&\wt{F}_l=\sum_{(m_1,m_2)\in\Z^2/\La_{d_1,d_2,d}}c_{d_2m_2-d_1m_1+l}
\Theta_{\La_{d_1,d_2,d},Q,mn>0;(-\frac{l}{d_0}+m_1,
\frac{l}{d_0}+m_2)}(x_0,x_3;\tau)=\\
&\sum_{(m,n)\in\Z^2}\eps(m,n)c_{d_2n-d_1m+l}\exp(\pi i\tau Q(m-\frac{l}{d_0},
n+\frac{l}{d_0})-
2\pi i [d_1x_1(m-\frac{l}{d_0})+d_2x_2(n+\frac{l}{d_0})])=\\
&\sum_{(m,n)\in\Z^2}\eps(-m,-n)c_{d_1m-d_2n+l}
\exp(\pi i\tau Q(m+\frac{l}{d_0},n-\frac{l}{d_0})+
2\pi i [d_1x_1(m+\frac{l}{d_0})+d_2x_2(n-\frac{l}{d_0})])
\end{align*}
where $\eps(m,n)=0$ unless $(m-\frac{l}{d_0}+\a(x_0))(n+\frac{l}{d_0}+\a(x_3))>0$ 
in which case $\eps(m,n)=\sign(m-\frac{l}{d_0}+\a(x_0))$.
Now an easy computation shows that the conditions (\ref{conds})
are equivalent to the system of equations
\begin{eqnarray}
\sum_{m\in\Z}c_{d_1m-d_2n_0+l}\exp(\pi i\tau Q(m+\frac{l}{d_0},n_0-\frac{l}{d_0})+
2\pi i d_1x_1(m+\frac{l}{d_0}))=0,\nonumber\\
\sum_{n\in\Z}c_{d_1m_0-d_2n+l}\exp(\pi i\tau Q(m_0+\frac{l}{d_0},n-\frac{l}{d_0})+
2\pi i d_2x_2(n-\frac{l}{d_0}))=0,\nonumber
\end{eqnarray}
where $m_0,n_0\in\Z$, $0\le l<d_0$.
Therefore, in the notation of Theorem \ref{linsys} we obtain
$$\wt{F}_l=\sum_{(m,n)\in\Z^2}\eps(-m,-n)c_{d_1m-d_2n+l}a_{m,n,l},$$
and the above system is equivalent to the condition \eqref{rowcol}.
This implies that one can replace the summation scheme defining $\wt{F}_l$ to the 
summation over $m,n\ge 0$ and $m,n<0$ (with signs ``minus" and ``plus", 
respectively). Hence, we obtain $\wt{F}_l=F_l$ (the latter series is defined
in the formulation of Theorem \ref{linsys}), i.e.,
\begin{equation}\label{MPFl-eq}
\lan e_{03}^*(0,l),MP(s_1,e,s_2)\ran=F_l.
\end{equation}

Now assume that the sections $s_1$ and $s_2$ have no common zeroes. 
Then we claim that formula \eqref{MPFl-eq} 
holds without any further genericity 
assumption on $x_1$ and $x_2$. Indeed, as we have seen in section 
\ref{masline-sec}, when the data $(\LL_0,\LL_1,\LL_2,s_1,s_2)$ vary in such a way 
that $s_1$ and $s_2$ have no common zeroes, the vector spaces $K_{s_1,s_2}$ can be 
viewed as fibers of a vector bundle on the space of parameters. Furthermore, the 
map $MP:K_{s_1,s_2}\to H^0(E,\LL_0)$ varies continuously with parameters.
Since $F_l$ is also a continuous function of $x_1$, $x_2$ and $(c_k)$ varying
in the vector bundle defined by \eqref{rowcol}, we derive that equation 
\eqref{MPFl-eq} holds whenever $s_1$ and $s_2$ have no common zeroes.

Finally, we are going to combine the result of section \ref{boundhom-sec} with 
Proposition \ref{mas} to derive the system of equations for $F_l$. We have
$$\delta_{-1}(e_{03}(0,l))=
\sum_{k\in\Z/d\Z}D_{k,l}e_{12}(0,k),$$
where
$$D_{k,l}=B_{s_1,s_2}(\th_{d_0\Z,l}(z,\frac{\tau}{d_0}), 
\th_{d\Z,-k}(z+x_0+x_1,\frac{\tau}{d})).$$
The divisor $Z(s_1)$ consists of the points $(-x_0+z(a),a\in\Z/d_1\Z)$
where
$$z(a)=\frac{\tau}{2}+\frac{1}{2d_1}+\frac{a}{d_1}$$
Therefore,
according to Lemma \ref{res} one has
\begin{align*}
&D_{k,l}=
\sum_{a\in \Z/d_1\Z}
\Res_{-x_0+z(a)}(\frac{2\pi i\th_{d_0\Z,l}(z,\frac{\tau}{d_0}) 
\th_{d\Z,-k}(z+x_0+x_1,\frac{\tau}{d})dz}
{\th_{d_1\Z}(z+x_0,\frac{\tau}{d_1})
\th_{d_2\Z}(z-x_3,\frac{\tau}{d_2})})=\\
&\sum_{a\in\Z/d_1\Z}\frac{2\pi i\th_{d_0\Z,l}(-x_0+z(a),\frac{\tau}{d_0})
\th_{d\Z,-k}(x_1+z(a),\frac{\tau}{d})}
{\th'_{d_1\Z}(z(a),\frac{\tau}{d_1})\th_{d_2\Z}(x_1+x_2+z(a),\frac{\tau}{d_2})}=\\
&\frac{2\pi i}{d_1\th'_{\Z}(\frac{d_1\tau+1}{2},d_1\tau)}
\sum_{a\in\Z/d_1\Z}
\frac{\th_{d_0\Z,l}(\frac{(d_2-d)x_1+d_2x_2}{d_0}+z(a),\frac{\tau}{d_0})
\th_{d\Z,-k}(x_1+z(a),\frac{\tau}{d})}
{\th_{d_2\Z}(x_1+x_2+z(a),\frac{\tau}{d_2})},
\end{align*}
where we denote
$\th'_{\Z}(z,\tau)=\frac{\partial}{\partial z}\th_{\Z}(z,\tau).$
It is easy to see that this formula is equivalent to the formula
for $D_{k,l}$ in Theorem \ref{linsys}.
Now Proposition \ref{mas} implies that
$$e=\delta_{-1}(\sum_l F_l e_{03}(0,l))=\sum_{k,l}D_{k,l}F_le_{12}(0,k).$$
Equating the coefficients with $e_{12}(0,k)$ we get the system of linear
equations on $F_l$.
It remains to notice that disjointness of the divisors of $s_1$ and $s_2$
is equivalent to the condition (\ref{disjoint}) of Theorem
\ref{linsys}.

\subsection{Examples of modular indefinite theta series}
\label{ex-mod}

Let us assume that $(d-d_i)|d$ for $i=1,2$. Choose an integer $f$ such that
$f|d$ and $(d-d_i)|f$ for $i=1,2$ and set 
$$c_k=\cases 1, \ k\equiv 0(f),\\ 0, \ k\not\equiv 0(f)\endcases$$
Then the conditions (\ref{conds}) boil down to the following pair of equations:
\begin{align*}
&\theta_{\frac{f}{d}\Z}(p_1x_1,p_1\tau)=0,\\
&\theta_{\frac{f}{d}\Z}(p_2x_2,p_2\tau)=0.  
\end{align*}
Hence, we can set
\begin{equation}\label{x1x2}
x_i=s_i\cdot\frac{f\tau}{2d}+r_i\frac{d-d_i}{2fd_i}, i=1,2,
\end{equation}
where $s_i$, $r_i$ are odd integers, to satisfy these equations.
Now Theorem \ref{linsys} implies that the series
$$-F_0=\sum_{(m,n)\in\Z^2, d_1m\equiv d_2n(f)}^{indef}
\exp(\pi i\tau Q(m,n)+2\pi i (d_1x_1m+d_2x_2n))$$
multiplied by the appropriate factor $\exp(\pi i c\tau)$ with $c\in\Q$, is
modular. More precisely, we can apply Theorem \ref{linsys} directly
only if the condition (\ref{disjoint}) is satisfied, i.e. the
corresponding section $s_1$ and $s_2$ have no common zeroes. Otherwise,
we first use Lemma \ref{divisor} to reduce to this case.

Let us denote by $f_0$ the least common multiple of $d-d_1$ and
$d-d_2$. Notice that the congruence
$d_1m-d_2n\equiv 0(f)$ implies that $(d-d_1)m$ and $(d-d_2)n$ belong to
$f_0\Z$. So we can take $\frac{(d-d_1)m}{f_0}$ and $\frac{(d-d_2)n}{f_0}$
as the new summation variables. Then we get
\begin{align*}
&-F_0=\sum_{(m,n)\in\Z^2, m\equiv n(\frac{f}{f_0})}^{indef}
\exp(\pi i\tau\frac{f_0^2}{d}(\frac{d_1}{d-d_1}m^2+
\frac{2d_1d_2}{(d-d_1)(d-d_2)}mn+\frac{d_2}{d-d_2}n^2)+\\
&\pi i\tau\frac{ff_0}{d}(\frac{d_1}{d-d_1}s_1m+\frac{d_2}{d-d_2}
s_2n)+\pi i\frac{f_0}{f}(r_1m+r_2n))
\end{align*}
To find the factor $\exp(\pi i c\tau)$ above recall that
for $N>0$ the functions 
$$\exp(\pi i\tau N\la^2)\th_{N\Z,i}(\la\tau+\mu,\frac{\tau}{N})$$
where $i\in\Z$, $\la,\mu\in\Q$, $\la>0$, and
$$\exp(\pi i\tau N/4)\th'_{\Z}(\frac{N\tau+1}{2},N\tau)$$
are modular. Hence,
$$\exp(\pi i \tau(d_0(\frac{f}{2dd_0}((d_2-d)s_1+d_2s_2)+\frac{1}{2})^2+
d(\frac{f}{2d}s_1+\frac{1}{2})^2-d_2(\frac{f}{2d}(s_1+s_2)+\frac{1}{2})^2
-\frac{d_1}{4}))D_{k,l}$$
are modular. Simplifying we conclude that
\begin{equation}\label{F0}
\exp(\pi i\tau\frac{f^2}{4d^2d_0}(d_1(d_2-d)s_1^2+2d_1d_2s_1s_2+d_2(d_1-
d)s_2^2))F_0
\end{equation}
is modular.

\noindent {\it Proof of Theorem \ref{modser}}. 
Let $a,b,c,p$ be positive integers such that $a|b$, $c|b$, 
$p|(b/a+1)$, $p|(c/a+1)$, $D=b^2-ac>0$. Then
we set $d_1=b(b+a)$, $d_2=b(b+c)$, $d=(b+a)(b+c)$, $d_0=D$,
$f=(b+a)(b+c)p/h$,
where $h$ is the greatest common divisor of $b/a+1$ and $b/c+1$,
so that we are in the situation considered above. 
We can rewrite the series (\ref{F0}) using the change of variables
$$\tau=\frac{ach^2}{b(b+a)(b+c)}\tau'.$$
Then the above argument implies the modularity of the series
$$q^{\frac{p^2ac(2bs_1s_2-as_1^2-cs_2^2)}{8D}}\cdot
\sum_{(m,n)\in\Z^2, m\equiv n(p)}^{indef}
\zeta_{2p}^{r_1m+r_2n}q^{bmn+a\frac{m^2+mps_1}{2}+
c\frac{n^2+nps_2}{2}},$$
where $\zeta_{2p}$ is a primitive root of $1$ of order $2p$ and
$r_1,r_2,s_1,s_2$ are odd integers. 
Note that we have
$$\zeta_{2p}^{r_1m+r_2n}=\cases\zeta_p^{\frac{r_1+r_2}{2}m}, \text{ if }
m\equiv n(2p),\\
-\zeta_p^{\frac{r_1+r_2}{2}m}, \text{ if } m\equiv n+p(2p).\endcases $$
Thus, the above series is equal to 
$$q^{\frac{p^2ac(2bs_1s_2-as_1^2-cs_2^2)}{8D}}\cdot
\sum_{(m,n)\in\Z^2, m\equiv n(p)}^{indef}(-1)^{\frac{n-m}{p}}
\zeta_{p}^{\frac{r_1+r_2}{2}m}q^{bmn+a\frac{m^2+mps_1}{2}+
c\frac{n^2+nps_2}{2}}.$$
Since $\frac{r_1+r_2}{2}$ can be an arbitrary integer we derive that
for any residue $r\in\Z/p\Z$ the series
$$q^{\frac{p^2ac(2bs_1s_2-as_1^2-cs_2^2)}{8D}}\cdot
\sum_{(m,n)\in\Z^2, m\equiv n\equiv r(p)}^{indef}
(-1)^{\frac{n-m}{p}}q^{bmn+a\frac{m^2+mps_1}{2}+
c\frac{n^2+nps_2}{2}}$$
is modular.
\ed

\noindent
{\it Proof of Theorem \ref{nonvan}}.
Let us denote the series $F_l$ considered above by $F_l(s_1,s_2)$
to show their dependence on a pair of odd integers $s_1,s_2$.
Note that using the change of variables from the proof of Theorem
\ref{modser} we can identify $f_{s_1,s_2}$ with $-F_0(s_1,s_2)$ multiplied
by some power of $q$.
By definition $F_l(s_1,s_2)=0$ unless $l\in(d-d_1)\Z+(d-d_2)\Z$.
On the other hand, an easy computation shows that for all
$l_1,l_2\in\Z$ one has
$$F_0(s_1+2\frac{h}{p}l_1,s_2+2\frac{h}{p}l_2)=\zeta\cdot
F_{(d-d_1)l_2-(d-d_2)l_1}(s_1,s_2)$$
where $\zeta$ is a root of unity.
The assumption in Theorem \ref{nonvan} is precisely the condition
(\ref{disjoint}) for the pair $x_1,x_2$ defined by (\ref{x1x2}).
By Theorem \ref{linsys}
at least one series among $F_l$ is non-zero which implies our
statement.
\ed

The simplest examples of identities that can be derived from the
above computations are obtained in the case $d_0=1$. Then $d=d_1+d_2-1$
so the conditions $(d-d_1)|d$, $(d-d_2)|d$ are satisfied only in the
following cases (assuming that $d_1\le d_2$): (i) $d_1=d_2=2$;
(ii) $d_1=2$, $d_2=3$; and (iii) $d_1=3$, $d_2=4$. In case (i) we obtain
formula \eqref{id1} taking $f=1$.
In case (ii) we get identities \eqref{id2}, \eqref{id3} and \eqref{id4}. 
More precisely, the case $f=4$ leads to \eqref{id2} while
in the case $f=2$ we get identities \eqref{id3} and \eqref{id4}
corresponding to the cases $s_1\equiv s_2(4)$ and $s_1\equiv -s_2(4)$.
In case (iii) above the assumption (\ref{disjoint}) of Theorem \ref{linsys}
is never satisfied, so we do not get any new identity.

One can also consider the degenerate case $d_0=d_1<d_2=d$ in
the above picture.
The coefficients of the Massey products in this case are given by the 
series \eqref{Jacform-ser}. Application of
the above analysis implies that
this is a meromorphic Jacobi form of weight $1$ as claimed in Theorem 
\ref{split-Jac-thm}.
On the other hand, we obtain
its expression in terms of theta functions. The case $a=1$ is well known
(see \cite{TM}, Section 486, or \cite{KP} (5.26),
or \cite{P-ap}).
In the case $a=2$ we get for any odd $s$
\begin{equation}\label{a2-eq}
\sum_{n\in\Z}\frac{q^{\frac{n^2+ns}{2}}}{1-q^{2n}u}=\varphi(q^2)^3
\frac{\sum_{n\in\Z} q^{2n^2+n(s-2)}u^n}{(\sum_{n\in\Z}(-1)^n
q^{n^2-n}u^n)(\sum_{n\in\Z}q^{2n^2+n(s-2)})},
\end{equation}
where $\varphi(q)=\prod_{n\ge 1}(1-q^n)$.

\subsection{String functions for $A_1^{(1)}$}\label{string} 

In the paper \cite{KP} the authors discovered a relation between
the string functions of irreducible highest weight representations
of affine Lie algebra $\gog$ of type $A_1^{(1)}$ and Hecke's 
modular forms for certain indefinite
quadratic modules. Namely, for every dominant
weight $\La$ and every weight $\la$ they define the {\it string function}
$c_{\la}^{\La}(\tau)$ which describes a 
part of the character of the irreducible
$\gog$-module with highest weight $\La$ (see \cite{KP}).   
One of the formulas they get for this function
(\cite{KP}, bottom of page 258)
is
$$\eta(\tau)^3c_{\la}^{\La}(\tau)=\sum_{k,n\in\frac{1}{2}\Z,
k\equiv n(\Z), k\ge|n|\text{or} -k>|n|}(-1)^{2k}\sign(k+\frac{1}{4})
q^{(m+2)(k+A)^2-m(n+B)^2},$$
where $m$ is the level of $\La$ (the value of $\La$ on the central
generator), $A$ and $B$ are rational numbers determined by $\La$ and
$\la$, such that $2(m+2)A\pm 2mB$ are odd integers.
Taking $k-n$ and $k+n$ as new summation variables we can rewrite this series
as follows:
$$q^{\frac{2bs_1s_2-s_1^2-s_2^2}{8(b^2-1)}}\cdot
\sum_{k,n\in\Z}^{indef}(-1)^{k+n}
q^{(m+1)kn+\frac{k^2+s_1k}{2}+\frac{n^2+s_2n}{2}},
$$
where $s_1=2(m+2)A+2mB$, $s_2=2(m+2)A-2mB$.
As shown in \cite{KP} this series is equal to
$$\Th^H_{\Z^2,Q;(A,B)},$$
where the quadratic form $Q$ on $\Z^2$ is given by $Q(x,y)=2(m+2)x^2-2my^2$.
Below we generalize this observation to other series considered in Theorem 
\ref{modser}.

\subsection{Proof of Theorem \ref{hecke-thm}}\label{hecke}

As we observed above one can replace the condition
$(m+\frac{1}{2})(n+\frac{1}{2})>0$ in the definition of the series
from Theorem \ref{modser} by any condition of the form
$(m+\a)(n+\b)>0$ where $\a,\b\not\in\Z$ (due to the vanishing
of the similar sum along the lines parallel to the generators on the cone).
So we can write this series as follows:
$$F=\sum_{(m,n)\in p\Z\oplus p\Z+\cc,(m+\frac{1}{4D})(n+\frac{1}{4D})>0}
\sign(m+\frac{1}{4D})(-1)^{m_0+n_0}q^{\frac{am^2+2bmn+cn^2}{2}},$$
where $\cc=(r+\frac{pc(bs_2-as_1)}{2D},r+\frac{pa(bs_1-cs_2)}{2D})$,
$(m_0,n_0)\in\Z^2$ is defined by $(m,n)=p(m_0,n_0)+\cc$.
Let us 
take $m'=m/p$ and $n'=(n+\frac{b}{c}m)/p$ as the new summation variables.
Note that if $(m,n)=p(m_0,n_0)+\cc$ then
$$(m',n')=(m_0,n_0+\frac{b}{c}m_0)+(\frac{r}{p}+\frac{acs}{2D},
r\frac{b/c+1}{p}+\frac{s_2}{2}),$$
where $s=-s_1+\frac{b}{a}s_2$. Thus, $(m',n')$ runs through
the coset $\Z^2+(\frac{r}{p}+\frac{acs}{2D},\frac{1}{2})$
intersected with the cone $(m'+\veps)(n'-\frac{b}{c}m'+\veps)>0,$
where $\veps>0$ is sufficiently small. Therefore, we can write
$$\pm F=\sum_{(m,n)\in\Z^2+(\frac{r}{p}+\frac{acs}{2D},\frac{1}{2}),
(m+\veps)(n-\frac{b}{c}m+\veps)>0}\sign(m+\veps)
(-1)^{n_0+(\frac{b}{c}+1)m_0}q^{p^2(cn^2-\frac{D}{c}m^2)},$$
where $(m,n)=(m_0,n_0)+(\frac{r}{p}+\frac{acs}{2D},\frac{1}{2})$
(the sign in front of $F$ depends on the parity of $r(b/c+1)/p+(s_2-1)/2$).
Now the lattice $\La$ defined in the formulation of the theorem
comes into play. Namely, splitting the above sum in two pieces
according to the parity of $n_0+(\frac{b}{c}+1)m_0$ we obtain
\begin{align*}
&\pm F=\sum_{(m,n)\in (\La+(\frac{r}{p}+\frac{acs}{2D},\frac{1}{2}))\cap
C}\sign(m+\veps) q^{\frac{p^2}{2}Q(m,n)}-\\
&\sum_{(m,n)\in (\La+(\frac{r}{p}+\frac{acs}{2D},-\frac{1}{2}))\cap
C}\sign(m+\veps) q^{\frac{p^2}{2}Q(m,n)},
\end{align*}
where $C$ is the cone $(m+\veps)(n-\frac{b}{c}m+\eps)>0$.
It is convenient to identify the lattice $\La$ with a $\Z$-submodule
in the field $K=\Q(\sqrt{D})$ by associating to $(m,n)\in\La$ 
the element $cn+m\sqrt{D}$. Then we have
$$\frac{1}{2}Q(m,n)=\frac{1}{c}\Nm(cn+m\sqrt{D}).$$
For two non-zero elements $k_1,k_2\in K$ let us denote by
$\lan k_1,k_2\ran=\Q_{>0}k_1+\Q_{>0}k_2$, $[k_1,k_2]=
\Q_{\ge 0}k_1+\Q_{\ge 0}k_2$, $\lan k_1,k_2]=\Q_{\ge 0}k_1+\Q_{>0}k_2$.
The intersection of a $\La$-coset with the
cone $C$ is equal to its intersection with the set
$[1,b+\sqrt{D}]\cup\lan -1,-b-\sqrt{D}\ran$.
Making the change of variables $(m,n)\mapsto (m,-n)$ in the second
sum of the above expression for $\pm F$ we can write:
\begin{align*}
&\pm F=\sum_{(m,n)\in
(\La+\cc)\cap([1,b+\sqrt{D}]\cup\lan -1,-b-\sqrt{D}\ran)}
\sign(m+\veps) q^{\frac{p^2}{2}Q(m,n)}-\\
&\sum_{(m,n)\in
(\La+\cc)\cap([-1,-b+\sqrt{D}]\cup\lan 1,b-\sqrt{D}\ran)}
\sign(m+\veps) q^{\frac{p^2}{2}Q(m,n)},
\end{align*}
where $\cc=(\frac{r}{p}+\frac{acs}{2D},\frac{1}{2})$. Since
$\La+\cc$ doesn't contain zero, the obtained series is equal to
$$\sum_{(m,n)\in
(\La+\cc)\cap(\lan b-\sqrt{D},b+\sqrt{D}]\cup
\lan -b-\sqrt{D},-b+\sqrt{D}])}
\sign(cn+m\sqrt{D}) q^{\frac{p^2}{2}Q(m,n)}.$$
Let us consider the totally positive element
$\eps=\frac{b+\sqrt{D}}{b-\sqrt{D}}\in K$.
Since $\Nm(\eps)=1$ the multiplication by $\eps$ preserves the
quadratic form $Q$ on $\La$. The direct computation shows that
the multiplication by $\eps$ preserves also $\La+\cc$, so we have  
$$\pm F=\sum_{(m,n)\in (\La+\cc)\cap(Q>0)/G_{\eps}}\sign(cn+m\sqrt{D})
q^{\frac{p^2}{2}Q(m,n)},$$
where $G_{\eps}\subset K^*$ is the subgroup generated by $\eps$.
Let $G$ be the subgroup of the group $\ker(\Nm:K^*\ra\Q^*)$
consisting of elements $k$ such that $k$ is totally positive and
$k(\LL+\cc)=\LL+\cc$. Then $G_{\eps}$ is a subgroup of finite
index in $G$, so we have
$$\pm F=|G/G_{\eps}|\cdot
\sum_{(m,n)\in (\La+\cc)\cap(Q>0)/G}\sign(cn+m\sqrt{D})
q^{\frac{p^2}{2}Q(m,n)}=|G/G_{\eps}|\cdot\Th^H_{\La,Q;\cc}(p^2\tau).$$
\ed

\subsection{Determinantal Jacobi forms}

We are going to compute explicitly the map $MP_{s_1,s_2}^{\det}$ (see section
\ref{det-Mas-sec}) in the setup of section \ref{masgen}. Note that 
we can trivialize the $1$-dimensional vector space $M$ using the canonical
bases in the relevant vector spaces.
First, let us compute the map
$$e_{s_1,s_2}:\sideset{}{^{d_0-1}}\We \Hom^1(L_1,L_2)^*\ra \Hom^1(L_1,L_2).$$
For every subset $S\subset\Z/d\Z$ of
cardinality $d-d_0+1$ let us denote by
$e_{12}^S$ the element in $\sideset{}{^{d_0-1}}\We \Hom^1(L_1,L_2)^*$
which is induced by the projection
$$\Hom^1(L_1,L_2)\ra \C^{\oplus d_0-1}:\sum_k y_k e_{12}(0,k)\mapsto
(y_k,k\not\in S).$$
Then
$$e_{s_1,s_2}(e_{12}^S)=\sum_{k\in S}c_{k,S}e_{12}(0,k)$$
where $(c_{k,S},k\in S)$
is the sequence of $(d-d_0)\times (d-d_0)$-minors (with signs)
of the
$(d-d_0)\times (d-d_0+1)$-matrix $R$ obtained by putting together
the $(d-d_1)\times (d-d_0+1)$-matrix
$R_1=(A_{ik}; i\in\Z/(d-d_1)\Z, k\in S)$
and the $(d-d_2)\times (d-d_0+1)$-matrix
$R_2=(B_{jk}; j\in\Z/(d-d_2)\Z, k\in S)$,
where $A_{ik}$ is zero unless $k\equiv i (d,d-d_1)$ in which case
$$A_{ik}=
\theta_{I_1,-\frac{k}{d}+n_1(k,i)}(p_1x_1,p_1\tau),$$
where $n_1(k,i)\in\Z/I_1$ is characterized by the congruence
$dn_1(k,i)\equiv k-i (d-d_1)$;
similarly $B_{jk}$ is zero unless $k\equiv j (d,d-d_2)$ in which case
$$B_{jk}=
\theta_{I_2,\frac{k}{d}+n_2(k,j)}(p_2x_2,p_2\tau),$$
where $n_2(k,j)\in\Z/I_2$ is characterized by the congruence
$dn_2(k,j)\equiv j-k (d-d_2)$.

Now we have
$$F_{l,S}=\lan e_{03}^*(0,l),MP_{s_1,s_2}^{\det}(e_{12}^S)\ran=
\sum_{k\in S} C_{l,k} c_{k,S},$$
where $C_{l,k}$ are defined by (\ref{Clk}).
In other words, $F_{l,S}$ 
is equal to the determinant of the
$(d-d_0+1)\times (d-d_0+1)$-matrix obtained by putting together
$R_1$, $R_2$ and
the row $(C_{l,k}, k\in S)$ of length $d-d_0+1$.

Let us recall the definition of Jacobi forms from \cite{GZ}. Let
$\La$ be a lattice, $\nn\cdot\nn'$ be a symmetric bilinear
form on $\La$ with values in $\Z$, $\Ga\subset\SL(2,\Z)$ be a
congruenz-subgroup. Let us denote $Q(\nn)=\nn\cdot\nn$ 
(this corresponds to $2Q$ in the notation of \cite{GZ}).
Then a meromorphic function
$f(\zz,\tau)$ on $\La_{\C}\times\HS$ is called a (meromorphic) Jacobi
form of weight $k$ with respect to $(\La,Q,\Ga)$
if the following equations hold:
$$f(\zz+\vv\tau+\ww,\tau)=(-1)^{\vv\cdot\ww}
\exp(-\pi i\tau Q(\vv)-2\pi i\vv\cdot\zz)f(\zz,\tau),$$
$$f(\frac{z}{c\tau+d},\frac{a\tau+b}{c\tau+d})=
(c\tau+d)^k\exp(\pi i\frac{c Q(\zz)}{c\tau+d})f(z,\tau),$$
for every $\vv,\ww\in\La$, $\left(\matrix a & b \\ c & d\endmatrix\right)
\in\Ga$. Similar to the case of modular forms one can extend the
above definition to half-integer weights $k$.

\begin{thm}\label{det-Jac-thm} 
The function $F_{l,S}$ is a Jacobi form of weight
$(d-d_0)/2+1$ with respect
to some sublattice $\La'\subset\La^+$, some congruenz-subgroup
$\Ga\subset\SL(2,\Z)$, and the quadratic form
$Q+Q_0$, where
$$\La^+=\{\nn=(n_0,n_1,n_2,n_3)\in\Z^4:\sum_i n_i=0, d_1n_1+(d_0-d_2)n_2+
d_0n_3=0\},$$
$$Q(\nn)=-d_1n_0n_1-d_2n_2n_3,$$
$$Q_0(\nn)=d(d_1n_1^2+d_2n_2^2).$$
\end{thm}

\Pf . We have the system of linear equations
$$\sum_{l} D_{k,l}F_{l,S}=c_{k,S}$$
determining $F_{l,S}$.
Using the functional equation for theta function we derive that
$D_{k,l}$ are Jacobi forms of weight $-1$ with respect
to some congruenz-subgroup of $\SL(2,\Z)$, some sublattice of $\La^+$
and the quadratic form $-Q$. On the other hand, $c_{k,S}$ are
Jacobi forms of weight $(d-d_0)/2$ with respect to the quadratic form
$Q_0$. Hence, $F_{l,S}$ are Jacobi forms of weight $(d-d_0)/2+1$
with respect to $Q+Q_0$.
\ed

\begin{rem} It is easy to see that the form $Q+Q_0$ can be written as follows:
$$(Q+Q_0)(\nn)=\frac{d_1d_2}{d_0}(x_1+x_2)^2+\frac{d(d_0-1)}{d_0}(d_1x_1^2+
d_2x_2^2).$$
In particular, it is always positive-definite for $d_0>1$.
In the case $d_0=1$ this form is degenerate which means that 
the corresponding functions $F_{l,S}$ have form $f(x_1+x_2,\tau)$,
where $f(z,\tau)$ is a Jacobi form.
\end{rem}

For example, for $d_1=d_2=2$, $d_0=1$ we have the unique choice of
$l$ and $S$.
Let us take $z_1=x_0,z_2=x_3$ as coordinates in $\La^+_{\C}$. Then
the corresponding function is equal (up to a sign)
to the following deteminant:
$$F(x_1,x_2;\tau)=\det\left(\matrix \Theta_{0}(x_1,x_2;\tau) &
\Theta_{-1}(x_1,x_2;\tau) &
\Theta_{1}(x_1,x_2;\tau) \\
\theta_0(6x_1,6\tau) & 
\theta_{-\frac{1}{3}}(6x_1,6\tau) & 
\theta_{\frac{1}{3}}(6x_1,6\tau) \\
\theta_0(6x_2,6\tau) & 
\theta_{\frac{1}{3}}(6x_2,6\tau) & 
\theta_{-\frac{1}{3}}(6x_2,6\tau) \endmatrix \right)
$$
where we denoted $\theta_r=\theta_{\Z,r}$,
\begin{align*}
&\Theta_c=\sum_{m-n\equiv c(3), (m+\a(2x_2-x_1))(n+\a(2x_1-x_2))>0}
\sign(m+\a(2x_2-x_1))\times\\
&\exp(\pi i\tau\frac{2}{3}(m^2+4mn+n^2)+4\pi i(mx_1+nx_2)).
\end{align*}
The function $F$
is a Jacobi form in $(x_1+x_2,\tau)$ of weight $2$ and index $2$.
Using the equation $D_{0,0}F=c_0$ and the addition formulas for
theta functions we derive the following identity:
$$F(x_1,x_2;\tau)=\theta_{\frac{1}{2}}(2(x_1+x_2)+\frac{1}{2},2\tau)^2
\cdot\frac{\eta^3(2\tau)\sum_{n\in\Z}\chi_3(n)q^{\frac{(4n+3)^2}{24}}}
{\theta(\frac{1}{2},2\tau)\theta_{\frac{1}{4}}(0,4\tau)}$$
where $\chi_3(n)$ is the non-trivial Dirichlet character modulo $3$.

\end{document}